\documentclass[11pt,twoside]{article}  

\usepackage[margin = 2.5cm]{geometry}

\usepackage{amsmath,amssymb,amsfonts,amsthm}
\usepackage{appendix}
\usepackage[dvipsnames]{xcolor}
\usepackage{siunitx}
\usepackage{graphicx}
\usepackage[numbers]{natbib}
\RequirePackage[colorlinks,citecolor=blue,urlcolor=blue]{hyperref}

\usepackage{enumitem}
\usepackage{subcaption}

\newtheorem{lem}{Lemma}
\newtheorem{thm}{Theorem}
\newtheorem{cor}{Corollary}

\newcommand{\survival}{\Psi}

\newcommand{\empsurv}{\hat{\Psi}}
\newcommand{\empsurvnull}{\hat{\Psi}_{0}}
\newcommand{\empsurvsig}{\hat{\Psi}_{1}}

\newcommand{\nulls}{\mathcal{H}_{0}}
\newcommand{\sigs}{\mathcal{H}_{1}}

\newcommand{\E}{\mathbb{E}}

\newcommand{\Var}{\mathrm{Var}}
\renewcommand{\P}{\mathbb{P}}
\newcommand{\R}{\mathbb{R}}

\newcommand{\N}{\mathbb{N}}
\newcommand{\FDR}{\mathrm{FDR}}

\newcommand{\FNR}{\mathrm{FNR}}
\newcommand{\Risk}{\mathcal{R}}

\newcommand{\FDP}{\mathrm{FDP}}
\newcommand{\FNP}{\mathrm{FNP}}
\newcommand{\FDPhat}{\widehat{\FDP}}
\newcommand{\one}{\boldsymbol{1}}
\newcommand{\discoveries}{\mathcal{I}}
\newcommand{\Xvals}{\mathcal{X}_{n}}
\newcommand{\poly}{\mathrm{poly}}
\newcommand{\event}{\mathcal{E}}
\newcommand{\eventwrap}[1]{\left\lbrace{#1}\right\rbrace}

\newcommand{\BH}{\mathrm{BH}}
\newcommand{\BC}{\mathrm{BC}}
\newcommand{\alg}{\mathrm{A}}
\newcommand{\iBH}{i_{\BH}}
\newcommand{\tBH}{t_{\BH}}

\newcommand{\tBC}{t_{\BC}}

\newcommand{\dgammaexpr}{D_{\gamma}}
\newcommand{\dgamma}[2]{\dgammaexpr\left({#1},{#2}\right)}

\newcommand{\rmax}{r_{\max}}

\newcommand{\widgraph}[2]{\includegraphics[keepaspectratio,width=#1]{#2}}

\newcommand{\mydefn}{\ensuremath{: \, =}}
\newcommand{\defn}{\mydefn}

\newcommand{\numobs}{\ensuremath{n}}
\newcommand{\kappan}{\ensuremath{\kappa_{\numobs}}}

\newcommand{\Zlower}{\ensuremath{Z_\ell}}
\newcommand{\Zupper}{\ensuremath{Z_u}}
\newcommand{\cfun}{\ensuremath{c}}

\newcommand{\card}[1]{\ensuremath{\mbox{card}(#1)}}
\newcommand{\Xseq}{\ensuremath{X_1^\numobs}}
\newcommand{\real}{\ensuremath{\mathbb{R}}}
\newcommand{\GG}{\ensuremath{\mbox{tGG}_\gamma}}
\newcommand{\newdisc}{\ensuremath{\discoveries}}
\newcommand{\RFUN}{\ensuremath{r_{\min}}}
\newcommand{\kappanstar}{\kappastar}
\newcommand{\betan}{\ensuremath{\beta_\numobs}}
\newcommand{\rn}{\ensuremath{r_\numobs}}
\newcommand{\kappastar}{\ensuremath{\kappa_*}}
\newcommand{\Event}{\ensuremath{\mathcal{D}}}
\newcommand{\Exs}{\ensuremath{\mathbb{E}}}
\newcommand{\Ftil}{\ensuremath{\widetilde{F}}}

\newcommand{\Bin}{\mathrm{Bin}}
\newcommand{\nulltxt}{\mathrm{null}}
\newcommand{\sigtxt}{\mathrm{signal}}

\newcommand{\nmin}{n_{\min}}
\newcommand{\nminl}{n_{\min,\ell}}
\newcommand{\nminu}{n_{\min,u}}
\newcommand{\rmin}{r_{\min}}
\newcommand{\taumin}{\tau_{\min}}
\newcommand{\tauminBC}{\tau_{\min,\BC}}
\newcommand{\tauminBH}{\tau_{\min,\BH}}
\newcommand{\qstar}{q_*}
\newcommand{\cstar}{c_*}

\long\def\comment#1{}

\makeatletter
\long\def\@makecaption#1#2{
        \vskip 0.8ex
        \setbox\@tempboxa\hbox{\small {\bf #1:} #2}
        \parindent 1.5em  
        \dimen0=\hsize
        \advance\dimen0 by -3em
        \ifdim \wd\@tempboxa >\dimen0
                \hbox to \hsize{
                        \parindent 0em
                        \hfil 
                        \parbox{\dimen0}{\def\baselinestretch{0.96}\small
                                {\bf #1.} #2
                                } 
                        \hfil}
        \else \hbox to \hsize{\hfil \box\@tempboxa \hfil}
        \fi
        }
\makeatother

\begin{document}

\begin{center}

{\bf{\LARGE{Optimal Rates and Tradeoffs in Multiple Testing}}}

\vspace*{.2in}

{\large{
\begin{tabular}{c}
Maxim Rabinovich$^\dagger$, Aaditya Ramdas$^{*,\dagger}$ \\ Michael
I. Jordan$^{*, \dagger}$, Martin J. Wainwright$^{*, \dagger}$ 
\\
\\
\texttt{\{rabinovich,aramdas,jordan,wainwrig\}@berkeley.edu}\\
\\
\end{tabular}
\begin{tabular}{ccc}
  Departments of Statistics$^*$ and EECS$^\dagger$, University of California, Berkeley
\end{tabular}
}}

\vspace*{.2in} \today
\vspace*{.2in}

\begin{abstract} 
Multiple hypothesis testing is a central topic in statistics, but
despite abundant work on the false discovery rate (FDR) and the
corresponding Type-II error concept known as the false non-discovery
rate (FNR), a fine-grained understanding of the fundamental limits of
multiple testing has not been developed. Our main contribution is to
derive a precise non-asymptotic tradeoff between FNR and FDR for a
variant of the generalized Gaussian sequence model.  Our analysis is
flexible enough to permit analyses of settings where the problem
parameters vary with the number of hypotheses $n$, including various
sparse and dense regimes (with $o(n)$ and $\mathcal{O}(n)$ signals).
Moreover, we prove that the Benjamini-Hochberg algorithm as well as
the Barber-Cand\`{e}s algorithm are both rate-optimal up to constants
across these regimes.
\end{abstract}
\end{center}


\section{Introduction}
\label{sec:intro}

The problem of multiple comparisons has been a central topic in
statistics ever since Tukey's influential 1953
book~\cite{tukey1953problem}.  In broad terms, suppose that one
observes a sequence of $\numobs$ independent random variables $X_1,
\ldots, X_\numobs$, of which some unknown subset are drawn from a null
distribution, corresponding to the absence of a signal or effect,
whereas the remainder are drawn from a non-null distribution,
corresponding to signals or effects. Within this framework, one can
pose three problems of increasing hardness: the \emph{detection}
problem of testing whether or not there is at least one signal; the
\emph{localization} problem of identifying the positions of the nulls
and signals; and the \emph{estimation} problem of returning estimates
of the means and/or distributions of the observations.  Note that
these problems form a hierarchy of difficulty: identifying the signals
implies that we know whether there is at least one of them, and
estimating each mean implies we know which are zero and which are
not. The focus of this paper is on the problem of localization.

There are a variety of ways of measuring type I errors for the
localization problem, including the \emph{family-wise error rate},
which is the probability of incorrectly rejecting at least one null,
and the \emph{false discovery rate} (FDR), which is the expected
ratio of incorrect rejections to total rejections.  An extensive
literature has developed around both of these metrics, resulting in
algorithms geared towards controlling one or the other.  Our focus
is the FDR metric, which has been widely studied, but for which
relatively little is known about the behavior of existing algorithms
in terms of the corresponding Type-II error concept, namely the \emph{false
  non-discovery rate} (FNR).\footnote{We
  follow~\citet{Ari16DistFreeFDR} in defining the FNR as the ratio of
  undiscovered to total non-nulls, which differs from the definition
  of Genovese and Wasserman~\cite{genovese2002operating}.}  Indeed, it
is only very recently that ~\citet{Ari16DistFreeFDR}, working within a
version of the sparse generalized Gaussian sequence model, established
asymptotic consistency for the $\FDR$-$\FNR$ localization problem.
Informally, in this framework, we receive $n$ independent observations
$X_1,\dots,X_n$, out of which $n - n^{1-\beta_n}$ are nulls, and the
remainder are non-nulls.  The null variables are drawn from a centered
distribution with tails decaying as
$\exp\big(-\frac{|x|^{\gamma}}{\gamma}\big)$, whereas the remaining
$n^{1-\beta_n}$ non-nulls are drawn from the same distribution shifted
by $(\gamma r_n \log n)^{1/\gamma}$.  Using this notation,
\citet{Ari16DistFreeFDR} considered the setting with fixed problem
parameters $r_n = r$ and $\beta_n = \beta$, and showed that when $r <
\beta < 1$, all procedures must have risk $\FDR + \FNR \to 1$. They
also showed that in the achievable regime $r > \beta > 0$, the
Benjamini-Hochberg (BH) is consistent, meaning that
$\FDR + \FNR \to 0$.  Finally, they proposed a new
``distribution-free'' method inspired by the knockoff procedure by
Barber and Cand\`{e}s~\cite{Bar15Knockoffs}, and they showed that the
resulting procedure is also consistent in the
achievable regime.

These existing consistency results are asymptotic.  To date there has been
no study of the important non-asymptotic questions that are of interest in
comparing procedures.  For instance, for a given FDR level, what is the best
possible achievable FNR?  What is the best-possible non-asymptotic
behavior of the risk $\FDR+\FNR$ attainable in finite samples?  And,
perhaps most importantly, non-asymptotic questions regarding whether
or not procedures such as BC and BH are \emph{rate-optimal} for the
FDR+FNR risk---remain unanswered. The main contributions of this paper
are to develop techniques for addressing such questions, and to essentially
resolve them in the context of the sparse generalized Gaussians model.

Specifically, we establish the tradeoff between FDR and FNR in finite
samples (and hence also asymptotically), and we use the tradeoff to
determine the best attainable rate for the $\FDR + \FNR$ risk.  Our
theory is sufficiently general to accommodate sequences of parameters
$(r_\numobs, \beta_\numobs)$, and thereby to reveal new phenomena that
arise when $r_n - \beta_n = o(1)$.  For a fixed pair of parameters
$(r, \beta)$ in the achievable regime $r > \beta$, our theory leads to
an explicit expression for the optimal rate at which FDR+FNR can
decay.  In particular, defining the $\gamma$-``distance''
$\dgamma{a}{b} \defn \big| a^{1/\gamma} - b^{1/\gamma} \big|^{\gamma}$
between pairs of positive numbers, we show that the equation
\begin{align*}
\kappa & = \dgamma{\beta + \kappa}{r}
\end{align*}
has a unique solution $\kappanstar$, and moreover that the combined
risk of any threshold-based multiple testing procedure $\newdisc$ is
lower bounded as $\Risk_\numobs(\newdisc) \gtrsim n^{-\kappanstar}$.
Moreover, by direct analysis, we are able to prove that both the
Benjamini-Hochberg (BH) and the Barber-Cand\`{e}s (BC) algorithms
attain this optimal rate.

At the core of our analysis is a simple comparison principle, and the
flexibility of the resulting proof strategy allows us to identify a
new critical regime in which $r_n - \beta_n = o(1)$, but the problem
is infeasible, meaning that if the $\FDR$ is driven to zero, then the
$\FNR$ must remain bounded away from zero.  Moreover, we are able to
study some challenging settings in which the fraction of signals is a
constant $\pi_1 \in (0,1)$ and not asymptotically vanishing, which
corresponds to the setting \mbox{$\beta_n = \frac{\log(1/\pi_1)}{\log
n}$,} so that $\beta_n \to 0$.  Perhaps surprisingly, even in
these regimes, the BH and BC algorithms continue to be optimal, though
the best rate can weaken from polynomial to subpolynomial in the
number of hypotheses $n$.


\subsection{Related work}
\label{subsec:related}

As noted above, our work provides a non-asymptotic generalization of
recent work by Arias-Castro and Chen~\cite{Ari16DistFreeFDR} on
asymptotic consistency in localization, using $\FDR+\FNR$ as the
notion of risk.  It should be noted that this notion of risk is
distinct from the asymptotic Bayes optimality under sparsity (ABOS)
studied in past work by~\citet{bogdan2011asymptotic} for Gaussian
sequences, and more recently by~\citet{Neu12FDRClass} for binary
classification with extreme class imbalance.  The ABOS results concern
a risk derived from the probability of incorrectly rejecting a single
null sample (false positive, or FP for short) and the probability of
incorrectly failing to reject a single non-null sample (false
negative, or FN for short).  Concretely, one has
$\Risk_{n}^{\mathrm{ABOS}} = w_{1} \cdot \mathrm{FP}~+ w_{2} \cdot
\mathrm{FN}$ for some pair of positive weights $(w_1, w_2)$ that need
not be equal.  As this risk is based on the error probability for a
single sample, it is much closer to misclassification risk or
single-testing risk than to the ratio-based $\FDR+\FNR$ risk studied
in this paper.

Using the notation of this paper, the work of~\citet{Neu12FDRClass}
can be understood as focusing on the particular setting $r = \beta$, a
regime referred to as the ``verge of detectability'' by these authors,
and with performance metric given by the Bayes classification risk,
rather than the combination of FDR and FNR studied here. In
comparison, our results provide additional insight into models that
are close to the verge of detectability, in that even when
$\beta_n=\beta$ is fixed, we can provide quantitative lower and upper
bounds on the FDR/FNR ratio as $r_\numobs \rightarrow \beta$ from above;
moreover, these bounds depend on how quickly $r_\numobs$ approaches
$\beta$.  These conclusions actually make it clear that a further
transition in rates occurs in the case where $r = \beta$ exactly for
all $n$, though we do not explore the latter case in depth. We suspect
that the methods developed in this paper may have sufficient precision
to answer the non-asymptotic minimaxity questions posed
by~\citet{Neu12FDRClass} as to whether any threshold-based procedure
can match the Bayes optimal classification error rate up to an
additive error $\ll \frac{1}{\log n}$.

The above line of work is complementary to the well-known asymptotic
results by Donoho and Jin~\citep{Don04HC,Don15HC} on phase transitions
in detectability using Tukey's higher-criticism statistic, employing
the standard type-I and type-II errors for testing of the single
global null hypothesis.  Note that Donoho and Jin use the generalized
Gaussian assumption directly on the PDFs, while our
assumption~\eqref{EqTailBounds} is on the survival function. Just as
in Arias-Castro and Chen~\cite{Ari16DistFreeFDR}, Donoho and Jin also
consider the asymptotic setting\footnote{We are not aware of any
  non-asymptotic results for detection akin to the results that the current paper
  provides for localization.} with $r_n = r$ and $\beta_n = \beta$,
which they sometimes call the ARW (asymptotic, rare and weak) model.

Our paper is also complementary to work on estimation, the most
notable result being the asymptotic minimax optimality of BH-derived
thresholding for denoising an approximately-sparse high-dimensional
vector~\citep{Abr06Adapt,Don06Asymptotic}. The relevance of our
results on the minimaxity of BH for approximately-sparse denoising
problems lies primarily in the use of deterministic thresholds as a
useful proxy for BH and other procedures that determine their
threshold in a manner that has complex dependence on the input
data~\citep{Don06Asymptotic}. Unlike the strategy
of~\citet{Don06Asymptotic}, which depends on establishing
concentration of the empirical threshold around the population-level
value, we use a more flexible comparison principle.  Deterministic
approximations to optimal FDR thresholds are also studied
by~\citet{Chi07Criticality} and~\citet{Gen06FDRWeights}. Other related
papers are discussed in Section~\ref{SecDisc}, when discussing
directions for future work.

The remainder of this paper is organized as follows. In
Section~\ref{SecBackground}, we provide background on the multiple
testing problem, as well as the particular model we consider.  In
Section~\ref{SecMain}, we provide an overview of our main results:
namely, optimal tradeoffs between FDR and FNR, which imply lower
bounds on the FDR+FNR risk, and optimality guarantees for the BH and
BC algorithms. In Section~\ref{SecProofs}, we prove our main results,
focusing first on the lower bounds and then using the ideas we have
developed to provide matching upper bounds for the well-known and
popular Benjamini-Hochberg (BH) procedure and the recent
Barber-Cand\`{e}s (BC) algorithm for multiple testing with FDR
control.  Proofs of some technical lemmas are given in the apepndices.


\section{Problem formulation}
\label{SecBackground}

 In this section, we provide background and a precise formulation of
 the problem under study.


\subsection{Multiple testing and false discovery rate}

Suppose that we observe a real-valued sequence $\Xseq \defn \{X_1,
\ldots, X_\numobs \}$ of $\numobs$ independent random variables.  
When the null hypothesis is true, $X_i$ is
assumed to have zero mean; otherwise, it is assumed that the mean of
$X_i$ is some unknown number $\mu_\numobs > 0$.  We introduce the
sequence of binary labels $\{H_1, \ldots, H_\numobs \}$ to encode
whether or not the null hypothesis holds for each observation; the
setting $H_i = 0$ indicates that the null hypothesis holds.  We define
\begin{align}
  \nulls & \defn \{ i \in [\numobs] \mid H_i = 0 \}, \quad \mbox{and}
  \quad \sigs \defn \{ i \in [\numobs] \mid H_i = 1 \}, \;
\end{align}
corresponding to the \emph{nulls} and \emph{signals}, respectively.
Our task is to identify a subset of indices that contains as many
signals as possible, while not containing too many nulls.

More formally, a testing rule $\discoveries: \real^\numobs \rightarrow
2^{[\numobs]}$ is a measurable mapping of the observation sequence
$\Xseq$ to a set $\discoveries(\Xseq) \subseteq [\numobs]$ of
\emph{discoveries}, where the subset $\discoveries(\Xseq)$ contains
those indices for which the procedure rejects the null hypothesis.
There is no single unique measure of performance for a testing rule
for the localization problem.  In this paper, we study the notion of
the \emph{false discovery rate} (FDR), paired with the
\emph{false non-discovery rate} (FNR).  These can be viewed
as generalizations of the type-I and type-II errors for single
hypothesis testing.

We begin by defining the false discovery proportion (FDP), and false
non-discovery proportion (FNP), respectively, as
\begin{align}
\FDP_\numobs(\discoveries) \defn \frac{\card{\discoveries(\Xseq) \cap
    \nulls}}{ \card{\discoveries(\Xseq)} \vee 1}, \quad \mbox{and}
\quad \FNP_\numobs(\discoveries) \defn \frac{\card{\discoveries(\Xseq)
    \cap \sigs}}{\card{\sigs}}.
\end{align}
Since the output $\discoveries(\Xseq)$ of the testing procedure is
random, both quantities are random variables.  The FDR and FNR are
given by taking the expectations of these random quantities---that is
\begin{align}
  \label{EqnFDRFNR}
\FDR_\numobs(\discoveries) & \defn \E\Big[\frac{ \card{
      \discoveries(\Xseq) \cap\nulls } }{\card {\discoveries(\Xseq)}
    \vee 1}\Big], \quad \mbox{and} \quad \FNR_\numobs(\discoveries)
\defn \E\Big[\frac{ \card{ \discoveries(\Xseq) \cap \sigs }
  }{\card{\sigs}} \Big],
\end{align}
where the expectation is taken over the random samples $\Xseq$. In
this paper, we measure the overall performance of a given procedure in
terms of its \emph{combined risk}
\begin{align}
  \label{EqnRisk}
\Risk_\numobs(\discoveries) \defn \FDR_\numobs(\discoveries) +
\FNR_\numobs(\discoveries).
\end{align}
Finally, when the testing rule $\discoveries$ under discussion is
clear from the context, we frequently omit explicit reference to this
dependence from all of these quantities.


\subsection{Tail generalized Gaussians model}

In this paper, we describe the distribution of the observations for
both nulls and non-nulls in terms of a \emph{tail generalized
  Gaussians model}. Our model is a variant of the generalized Gaussian
sequence model studied in past work~\citep{Ari16DistFreeFDR,Don04HC};
the only difference is that whereas a $\gamma$-generalized Gaussian
has a density proportional to
$\exp\big(-\frac{|x|^{\gamma}}{\gamma}\big)$, we focus on
distributions whose tails are proportional to
$\exp\big(-\frac{|x|^{\gamma}}{\gamma}\big)$. This alteration is in
line with the asymptotically generalized Gaussian (AGG) distributions
studied by~\citet{Ari16DistFreeFDR}, with the important caveat that
our assumptions are imposed in a non-asymptotic fashion.

For a given degree $\gamma \geq 1$, a $\gamma$-tail generalized
Gaussian random variable with mean 0, written as $G\sim \GG(0)$, has a
survival function \mbox{$\survival(t) \defn \P \big( G \geq t \big )$}
that satisfies the bounds
\begin{align}\label{EqTailBounds}
\frac{e^{\frac{-|t|^\gamma}{\gamma}}}{\Zlower} ~\leq~
\min\{\survival(t), 1-\survival(t)\} ~\leq~
\frac{e^{\frac{-|t|^\gamma}{\gamma}}}{\Zupper},
\end{align}
for some constants $\Zlower > \Zupper > 0$.  (Note that $t \mapsto
\Psi(t)$ is a decreasing function, and becomes smaller than
$1-\Psi(t)$ at the origin.)  As a concrete example, a $\gamma$-tail
generalized Gaussian with $\Zlower = \Zupper = 1$ can be generated by
sampling a standard exponential random variable $E$ and a Rademacher
random variable $\varepsilon$ and putting $G = \varepsilon \big(\gamma
E\big)^{1/\gamma}$.  We use the terminology ``tail generalized Gaussian''
because of the following connection: the survival function of a $2$-tail
Gaussian random variable is on the order of $\exp(-|x|^{2}/2)$,
whereas that of a Gaussian is on the order of $\frac{1}{\poly(x)}
\exp(-x^2/2)$.  In particular, this observation implies a
$\mathrm{tGG}_{2}$ random variable has tails that are equivalent to a
Gaussian in terms of their exponential decay rates.

In terms of this notation, we
assume that each observation $X_i$ is distributed as
\begin{align}
X_i & \sim \begin{cases} \GG(0) & \mbox{if $i \in \nulls$}
  \\ \GG(0) + \mu_\numobs & \mbox{if $i \in \sigs$,}
  \end{cases}
\end{align}
where our notation reflects the fact that the mean shift $\mu_\numobs$
is permitted to vary with the number of observations $\numobs$.  See
Section~\ref{SecScaling} for further discussion of the scaling of the
mean shift.

\subsection{Threshold-based procedures}
\label{subsec:thresh-procedures}

Following prior work~\citep{Ari16DistFreeFDR,Don04HC}, we restrict
attention to testing procedures of the form
\begin{align}\label{EqThreshBased}
\discoveries(\Xseq) = \big \{ i \in [\numobs] \mid X_i \geq T_n(\Xseq)
\big \},
\end{align}
where $T_n(\Xseq) \in \real_+$ is a data-dependent threshold.  We
refer to such methods as \emph{threshold-based procedures}. The BH and BC procedures both belong to this class. 
Moreover, from an intuitive standpoint, the observations are exchangeable in the absence of prior information, and we are considering testing between a single unimodal null distribution and a single positive shift of that distribution. 
In this setting, it is hard to imagine that an optimal procedure would ever reject the hypothesis corresponding to one observation while rejecting a hypothesis with a smaller observation
value. Threshold-based procedures therefore appear to be a very reasonable class to focus on. 

It will be convenient to reason about the performance
metrics associated with rules of the form
\begin{align}
\label{eq:discoveries-t-def}
\newdisc_t(\Xseq) & = \big \{ i \in [\numobs] \mid X_i \geq t \big \},
\end{align}
where $t > 0$ is a pre-specified (fixed, non-random) threshold.  In this case, we adopt
the notation $\FDR_\numobs(t)$, $\FNR_\numobs(t)$ and
$\Risk_\numobs(t)$ to denote the metrics associated with
the rule $\Xseq \mapsto \newdisc_t(\Xseq)$.

\subsection{Benjamini-Hochberg (BH) and Barber-Cand\`{e}s (BC) procedures} 

Arguably the most popular threshold-based procedure that provably
controls FDR at a user-specified level $q_{n}$ is the
\emph{Benjamini-Hochberg} (BH) procedure. More recently,
~\citet{Ari16DistFreeFDR} proposed a method that we refer to as the
\emph{Barber-Cand\`{e}s} (BC) procedure.  Both algorithms are based on
estimating the $\FDP_\numobs$ that would be incurred at a range of
possible thresholds and choosing one that is as large as possible
(maximizing discoveries) while satisfying an upper bound linked to
$q_\numobs$ (controlling $\FDR_\numobs$).  Further, they both only
consider thresholds that coincide with one of the values $\Xseq$,
which we denote as a set by $\Xvals = \big \lbrace X_{1}, \dots,
X_\numobs \big \rbrace$. The data-dependent threshold for both can be
written as
\begin{align}
\label{eq:alg-thresh-def}
t_\numobs \big ( X_{1},\dots,X_\numobs \big ) & = \max \big\lbrace t
\in \Xvals \colon \FDPhat_\numobs\big ( t \big ) \leq
q_\numobs\big\rbrace .
\end{align}
The two algorithms differ in the estimator $\FDPhat_\numobs\big ( t
\big ) $ they use. The BH procedure assumes access to the true null
distribution through its survival function $\survival$ and sets
\begin{align}
  \label{eq:bh-fdp-def}
\FDPhat_\numobs^{\BH}\big ( t \big ) = \frac{\survival\big ( t \big )
}{\#\big ( X_{i} \geq t \big ) / n},~ \mbox{ for } t\in\Xvals .
\end{align}
The BC procedure instead estimates the survival function
$\survival(t)$ from the data and therefore does not even need to know
the null distribution. This approach is viable when $\#\big ( X_{i}
\leq -t \big ) /n$ is a good proxy for $\survival\big ( t \big ) $,
which our upper and lower tail bounds guarantee; more typically, the
BC procedure is applicable when the null distribution is (nearly)
symmetric, and the signals are shifted by a positive amount (as they
are in our case).  Then, the BC estimator is given by
\begin{align}
  \label{eq:bc-fdp-def}
\FDPhat_\numobs^{\BC}\big ( t \big ) = \frac{ \big[ \#\big ( X_{i}
    \leq -t \big ) + 1\big]/n}{\#\big ( X_{i} \geq t \big ) / n},
\quad \mbox{for $t \in \Xvals$.}
\end{align}
With these definitions in place, are now ready to describe our main results.


\section{Main results}
\label{SecMain}

We now turn to a statement of our main results, along with some
illustrations of their consequences.  Our first main result
(Theorem~\ref{ThmLower}) characterizes the optimal tradeoff between FDR and FNR for any testing procedure. 
By optimizing this tradeoff, we obtain a lower bound on the combined FDR and FNR of any testing procedure (Corollary~\ref{CorCool}).
Our second main result (Theorem~\ref{ThmUpper}), shows that BH achieves the optimal FDR-FNR tradeoff up to constants and that
BC almost achieves it. In particular, our result implies that with the proper choice of target FDR, both BH and BC can achieve the optimal combined FDR-FNR rate (Corollary~\ref{CorCoolUpper}).


\subsection{Scaling of sparsity and mean shifts}
\label{SecScaling}

We study a sparse instance of the multiple testing
problem in which the number of signals is assumed to be small relative
to the total number of hypotheses.  In particular, motivated by related work
in multiple hypothesis testing~\citep{Ari16DistFreeFDR,Don04HC,Don15HC,Jin14RareWeak},
we assume that the number of signals scales as
\begin{align}
  \label{EqnSparseSignal}
  \card{\sigs} & = m_\numobs = n^{1 - \beta_\numobs} \quad \mbox{for some
    $\beta_\numobs \in (0,1)$.}
\end{align}
Note that to the best of our knowledge, all previous results in the literature
assume that $\beta_\numobs = \beta$ is actually independent of
$\numobs$.  In this case, the sparsity
assumption~\eqref{EqnSparseSignal} implies that all but a polynomially
vanishing fraction of the hypotheses are null.  In contrast, as
indicated by our choice of notation, the set-up in this paper allows
for a sequence of parameters $\betan$ that can vary with the number
of hypotheses $\numobs$.  In this way, our framework is flexible
enough to handle relatively dense regimes (e.g., those with
$\frac{\numobs}{\log \numobs}$ or even $\mathcal{O}(n)$ signals).

The non-null hypotheses are distinguished by a positively shifted mean
$\mu_\numobs > 0$.  It is natural to parameterize this mean shift in
terms of a quantity $r_\numobs > 0$ via the relation
\begin{align}
  \label{EqnMeanShift}
\mu_\numobs = \big ( \gamma r_\numobs \log{n} \big )^{1/\gamma}.
\end{align}
As shown by Arias-Castro and Chen~\citep{Ari16DistFreeFDR}, when the
pair $(\beta, r)$ are fixed such that $r < \beta$, the problem is
asymptotically infeasible, meaning that there is no procedure such
that $\Risk_\numobs(\discoveries) \rightarrow 0$ as $\numobs
\rightarrow \infty$.  Accordingly, we focus on sequences
$(\beta_\numobs, r_\numobs)$ for which $r_\numobs > \beta_\numobs$.
Further, even though the asymptotic consistency boundary of $r <
\beta$ versus $r > \beta$ is apparently independent of $\gamma$, we
will see that the rate at which the risk decays to zero is determined
jointly by $r, \beta$ and $\gamma$.


\subsection{Lower bound on any threshold-based procedure}

In this section, we assume :
\begin{subequations}
\label{EqnAnnoying}
\begin{align}
\label{EqnAnnoying1}
\beta_\numobs ~\stackrel{(i)}{\geq}~ \frac{\log 2}{\log n} &\quad \iff \quad \numobs^{1-\beta_n} ~\leq~ \numobs/2, 
\quad \text{ and }\\
\label{EqnAnnoying2}
\max\{\betan,\frac1{\log^{\frac{\gamma - 1/2}{\gamma}} \numobs}\} &\quad \stackrel{(ii)}{<} \quad r_\numobs  \quad \stackrel{(iii)}{<}  \quad \rmax  \quad \text{ for some constant $\rmax <1$}.
\end{align}
\end{subequations}
Condition (i) requires that the proportion $\pi_1$ of non-nulls is at most 1/2. 
Condition (ii) asserts that the natural requirement of $r_n > \beta_n$ is not enough,
but further insists that  $r_n$ cannot approach zero too fast. 
The constants $\log{2}$ and $\frac{\gamma - 1/2}{\gamma}$ are
somewhat arbitrary and can be replaced, respectively, by $\log
\frac1{\pi_{\max}}$ for any $0 < \pi_{\max} < 1$ and $\frac{\gamma -
  1 + \rho}{\gamma} $ for any $\rho > 0$, but we fix their values in
order not to introduce unnecessary extra parameters.
As for condition (iii), although the assumption $r_\numobs
< 1$ is imposed because the problem becomes qualitatively
easy for $r_\numobs \geq 1$, the assumption that it is bounded away from one
 is a technical convenience that simplifies some of our proofs. 

Our analysis shows that the FNR behaves differently depending on the
closeness of the parameter $r_\numobs$ to the boundary of feasiblity
given by $\beta_\numobs$.  In order to characterize this closeness, we
define
\begin{align}
\label{eq:rn-minus}
\RFUN = \RFUN(\kappa_n) & \defn \begin{cases} \beta_\numobs + \kappa_{\numobs}
  + \frac{\log\frac{1}{6\Zlower}}{\log{n}} &~\text{if}~
  \kappan \leq 1 - \beta_n - \frac{\log \frac3{\log{16}}}{\log
    n}, \\
  1 + \frac{\log\frac{1}{24\Zlower}}{\log{n}}
  &~\text{otherwise.}
 \end{cases} 
\end{align}
Here $\kappan$ is to be interpreted as the ``exponent'' of a target
FDR rate $q_{n}$, in the sense that $q_{n} = n^{-\kappa_{n}}$. The
rate $q_{n}$ may differ from the actual achieved $\FDR_{n}$, but it is
nonetheless useful for parameterizing the quantities that enter into
our analysis. When we need to move between $q_{n}$ and $\kappa_{n}$,
we shall write $\kappan = \kappan(q_{n}) = \frac{\log(1/q_n)}{\log n}$
and $q_{n} = q_{n}(\kappa_{n}) = n^{-\kappan}$.  For mathematical
convenience, we wish to have the target FDR $q_{n}$ to be bounded away
from one, and we therefore impose one further technical but
inessential assumption in this section:
\begin{align}\label{EqnAnnoying3}
q_\numobs \leq \min\big\lbrace \frac{1}{24},~\frac{1}{6\Zlower}\big\rbrace \quad \iff \quad \kappan \geq \frac{\log \max\big\lbrace 24,~6\Zlower\big\rbrace}{\log{n}}.
\end{align}

\begin{figure}[t]
  \begin{center}
\begin{tabular}{ccc}
  \widgraph{0.45\textwidth}{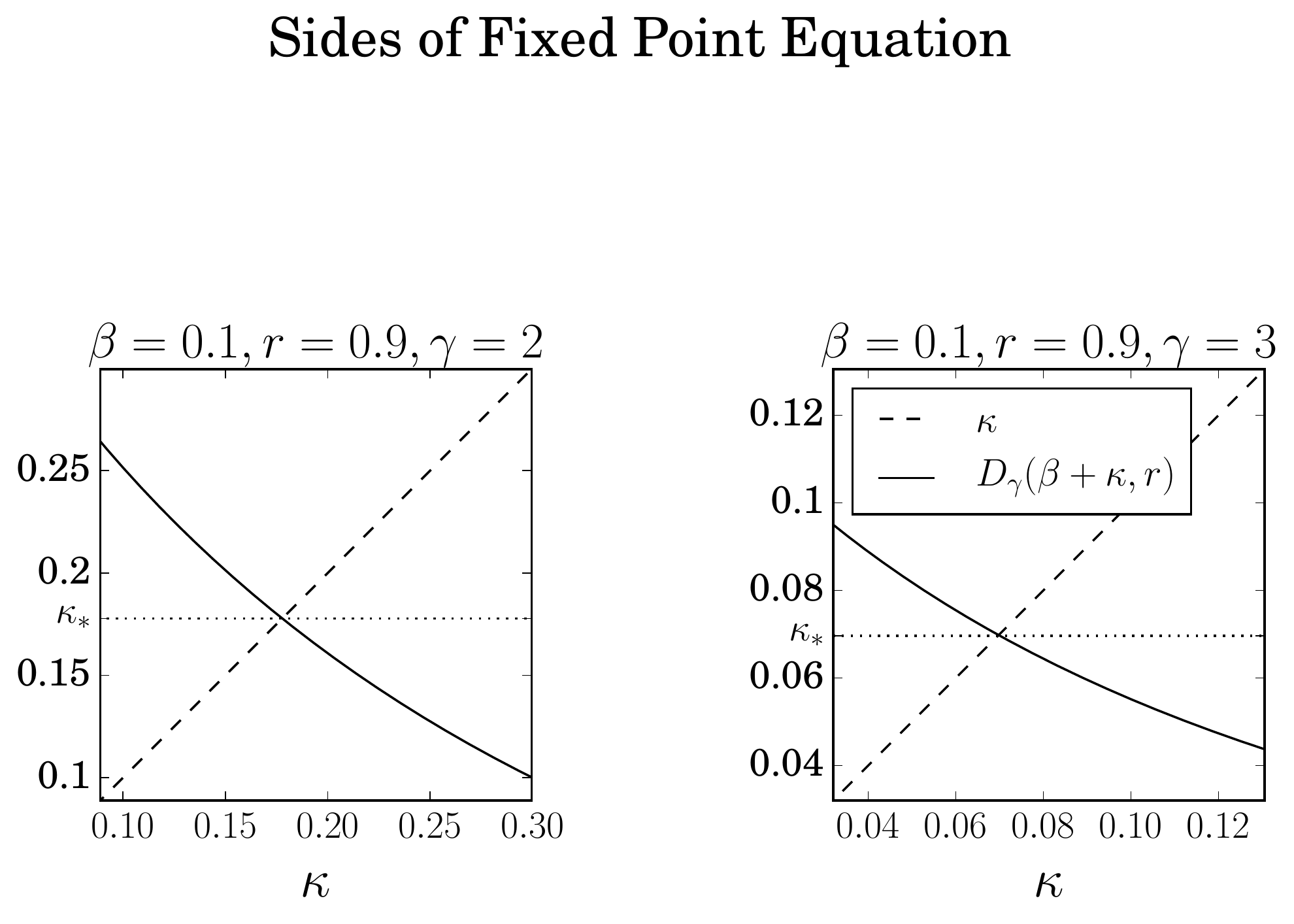} & &
  \widgraph{0.45\textwidth} {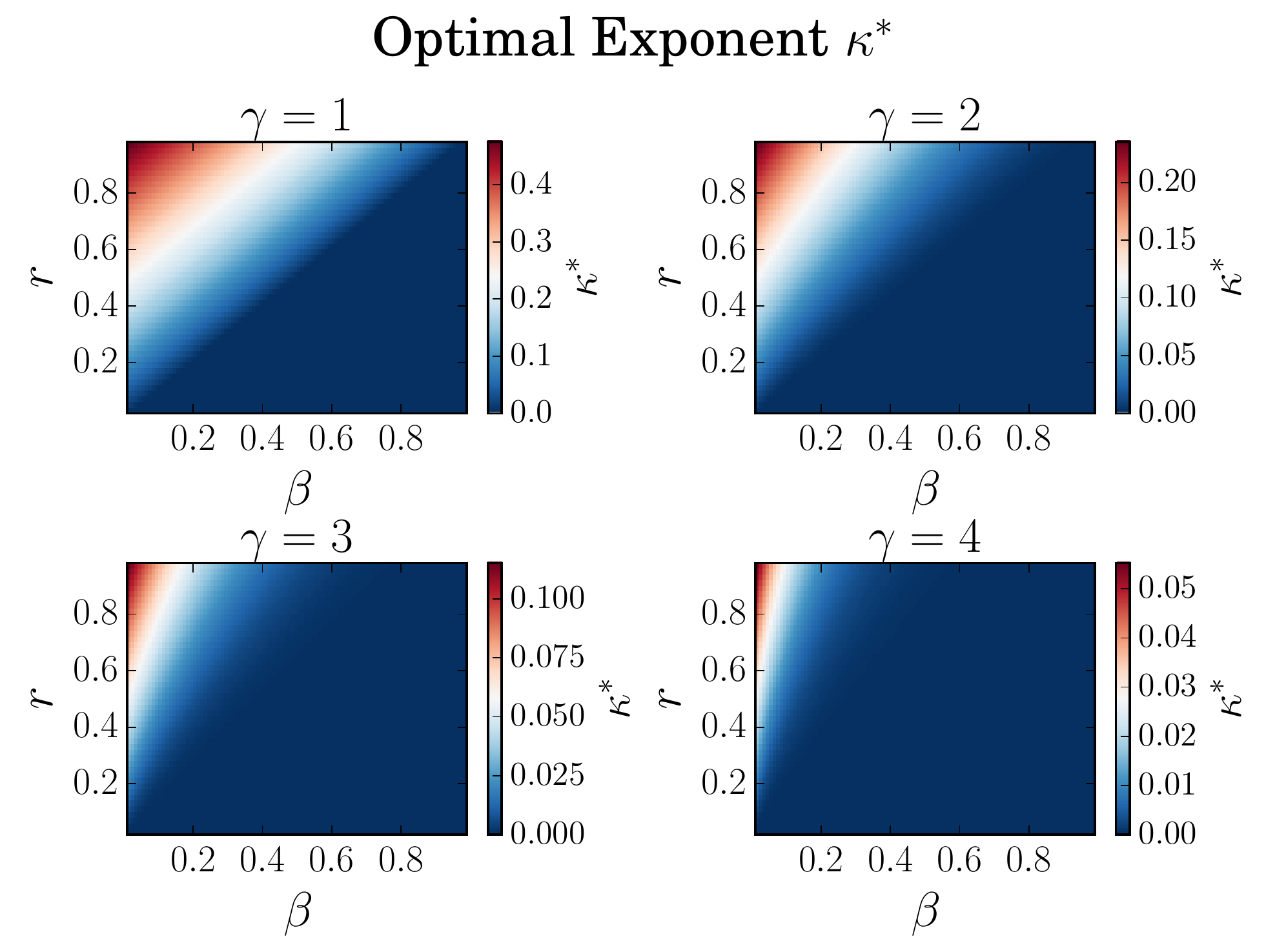} \\
  (a) & & (b)
\end{tabular}
\caption{ Visualizations of the fixed-point equation~\eqref{EqnCool}.
  (a) Plots comparing the left- and right-hand sides of the
  fixed-point equation.  (b) The optimal exponent $\kappa^{\ast}$ as a
  function of $r$ and $\beta$.}
\label{FigKappaStar}
  \end{center}
\end{figure}

\noindent The theorem that follows will apply to all sample sizes $n >
\nminl$ (subscript $\ell$ for lower), where
\begin{align}
\nminl ~&:=~ \min\left\{ n \in \N : \exp \left(-\frac{n^{1 -
    \rmax}}{24(\Zlower \vee 1)} \right) \leq \frac14 \right\} \; = \;
\left\lfloor [24(\Zlower \vee 1) \log 4]^{\frac1{1-\rmax}}
\right\rfloor. \label{EqNminLower}
\end{align}
Finally, for $\gamma \in [1, \infty)$ and non-negative numbers $a, b >
  0$, let us define the associated $\gamma$-``distance'':
\begin{align}
\label{eq:dgamma-defn}
\dgamma{a}{b} & \defn \big| a^{1/\gamma} - b^{1/\gamma}
\big|^{\gamma}.
\end{align}
Our first main theorem states that for $\rn > \RFUN(\kappa_n)$, the $\FNR$ decays 
as a power of $1/\numobs$, with exponent specified by the $\gamma$-distance.
\begin{thm}
\label{ThmLower}
Consider the $\gamma$-tail generalized Gaussians testing
problem with sparsity $\betan$ and signal level $\rn$ satisfying
conditions~\eqref{EqnAnnoying1}, and~\eqref{EqnAnnoying2}, and
 with sample size $\numobs > \nminl$ from
definition~\eqref{EqNminLower}.  Then, for any choice of exponent $\kappan
\in (0,1)$ satisfying condition~\eqref{EqnAnnoying3}, there exists a minimum signal strength $\RFUN(\kappan)$ from definition~\eqref{eq:rn-minus}, 
such that any threshold-based procedure $\newdisc$
that satisfies \mbox{$\FDR_\numobs(\newdisc) \leq \numobs^{-\kappan}$} must have its $\FNR$ lower bounded as
\begin{align}
  \label{EqnLower}
  \FNR_\numobs(\newdisc) & \geq \begin{cases} \frac{1}{32} & \mbox{if
      $\rn \in \big[\betan, \RFUN \big]$} \\
\cfun(\betan, \gamma) \; n^{-\dgamma{\betan + \kappan}{\rn}} &
\mbox{otherwise,}
  \end{cases}
\end{align}
where $\cfun(\betan, \gamma) \defn c_0 \exp \big(c_1 \betan^{\frac{1 -
    \gamma}{\gamma}} \big)$, with $(c_0, c_1)$ being positive constants
depending only on $(\Zlower, \Zupper, \gamma)$.
\end{thm}
\noindent The proof of this theorem is provided in Section~\ref{SecProofThmLower}. 
Note that the theorem holds for any choice of $\kappan \in (0,1)$.  In
the special case of constant pairs $(\beta, r)$, this choice can be
optimized to achieve the best possible lower bound on the
 risk $\Risk_\numobs(\newdisc) = \FDR_\numobs(\newdisc) +
\FNR_\numobs(\newdisc)$, as summarized below.  
%
\begin{cor}
\label{CorCool}
When $r > \beta$, let $\kappanstar = \kappanstar(\beta, r, \gamma) >
0$ be the unique solution to the equation
\begin{align}
  \label{EqnCool}
\kappa & = \dgamma{\beta + \kappa}{r}.
\end{align}
Then the combined risk of any threshold-based multiple testing
procedure $\newdisc$ is lower bounded as
\begin{align}
  \label{EqnCleanLower}
\Risk_\numobs(\newdisc) & \gtrsim n^{-\kappanstar},
\end{align}
where $\gtrsim$ denotes inequality up to a pre-factor independent of
$\numobs$.
\end{cor}
\noindent  The proof of this corollary is provided in Section~\ref{SecProofCorCool}. 
Figure~\ref{FigKappaStar} provides an illustration of the predictions
in Corollary~\ref{CorCool}.  In particular, panel (a) shows how the
unique solution $\kappastar$ to equation~\eqref{EqnCool} is determined
for varying settings of the triple $(r, \beta, \gamma)$.  Panel (b)
shows how $\kappanstar$ varies over the interval $(0, 0.5)$, again for
different settings of the triple $(r, \beta, \gamma)$.  As would be
expected, the fixed point $\kappanstar$ increases as a function of the
difference $r - \beta > 0$.

\subsection{Upper bounds for some specific procedures}

Thus far, we have provided general lower bounds applicable to any
threshold procedure.  We now turn to the complementary question---how
do these lower bounds compare to the results achievable by the
BH and BC algorithms introduced in
Section~\ref{subsec:thresh-procedures}?  Remarkably, we find that up
to the constants defining the prefactor, both the BH and BC procedures
achieve the minimax lower bound of Theorem~\ref{ThmLower}.

We state these achievable results in terms of the fixed point
$\kappanstar$ from equation~\eqref{EqnCool}.  Moreover, they apply to
 all problems with sample size $n > \nminu$ (subscript $u$ for upper), where
\begin{align}
\nminu ~&:=~ 
\min \left\{ n \in \N :  \exp \left(-\frac{n^{1 - \rmax}}{24}\right) 
\leq \frac{1}{\Zupper  n} \right\} \nonumber \\
&=~ 
\min \left\{ n \in \N :  n
~\geq~ [24\log(\Zupper  n)]^{\frac1{1-\rmax}} \right\} \label{EqNminUpper}.
\end{align}

In order to state our results cleanly, let us introduce the constants
\begin{align}
\label{EqZetaDef}
c_{\BH} ~:=~ \frac{\Zupper}{36\Zlower}, \quad c_{\BC} ~:=~
\frac{\Zupper}{48\Zlower}, \quad \mbox{and} \quad \zeta ~:=~
\max\big\lbrace 6\Zlower,~\frac{1}{6\Zlower}\big\rbrace,
\end{align}
and require in particular that $r_{\numobs} \geq
\rmin\left(\kappan\left(c_{\mathrm{A}}q_{\numobs}\right)\right)$ for
algorithm $\mathrm{A} \in \{\BH,\BC\}$.  Note that $c_{\mathrm{A}} <
1$ since $Z_\ell \geq Z_u$ by definition, and that the introduction of
$c_{\mathrm{A}}$ into the argument of $\rmin$ only changes the minimum
allowed value of $r_{n}$ by a conceptually negligible amount of
$\mathcal{O}\big(\frac{1}{\log n}\big)$.

Lastly, we note that BC requires an additional mild condition that the
number of non-nulls $n^{1-\beta_n}$ is large relative to the target
FDR $q_n = n^{-\kappan}$ (otherwise, in some sense, the problem is too
hard if there are too few non-nulls and a very strict target
FDR). Specifically, we need that both quantities cannot simultaneously
be too small, formalized by the assumption:
\begin{align}
  \label{eq:qn-decay}
\exists n_{\min,\BC}~~\text{such that for all}~~ n \geq n_{\min,\BC}
~~\text{we have}~~ \frac{3c_{\BC}}{4}\cdot \frac{q_\numobs}{\log
  \frac{1}{q_\numobs}} \cdot n^{1-\beta_n} ~\geq~ 1.
\end{align}
We note that when $r_{n} = r$ and $\beta_n = \beta$ are constants,
this decay condition is satisfied by $q_\numobs =
\numobs^{-\kappastar}$.

Our second main theorem delivers an optimality result for the BH and
BC procedures, showing that under some regularity conditions, their
performance achieves the lower bounds in Theorem~\ref{ThmLower} up to
constant factors.
\begin{thm}
\label{ThmUpper}
Consider the $\betan$-sparse $\gamma$-tail generalized Gaussians
testing problem with target FDR level $q_n$ upper bounded as in
condition~\eqref{EqnAnnoying3}.
\begin{enumerate}
\item[(a)] Guarantee for BH procedure: Given a signal strength
  $r_{\numobs} \geq \rmin (\kappan (c_{\BH} q_{\numobs}))$
  and sample size $\numobs > \nminu$ as in condition~\eqref{EqNminUpper},
  the BH procedure satisfies the bounds
  \begin{align}
    \label{EqnBHBound}
\FDR_\numobs \leq q_\numobs ~~ \text{and} ~~ \FNR_\numobs \leq
\frac{2\zeta_{\BH}^{2\beta_\numobs^{\frac{1 -
        \gamma}{\gamma}}}}{\Zupper} \cdot n^{-\dgamma{\beta_\numobs +
    \kappan}{r_\numobs}}, \qquad \mbox{where $\zeta_\BH :=
  \frac{\zeta}{c_\BH}$.}
  \end{align}
  \item[(b)] Guarantee for BC procedure: Given a signal strength
    $r_{\numobs} \geq \rmin (\kappan (c_{\BC} q_{\numobs}))$ and
    sample size $\numobs > \max\{n_{\min,\BC}, \nminu\}$ as in
    condition~\eqref{eq:qn-decay}, the BC procedure satisfies the
    bounds
  \begin{align}
    \label{EqnBCBound}
\FDR_\numobs \leq q_\numobs ~~ \text{and} ~~ \FNR_\numobs \leq
\frac{2\zeta_{\BC}^{2\beta_\numobs^{\frac{1 -
        \gamma}{\gamma}}}}{\Zupper} \cdot n^{-\dgamma{\beta_\numobs +
    \kappan}{r_\numobs}} + q_n, \qquad \mbox{where $\zeta_\BC :=
  \frac{\zeta}{c_\BC}$.}
  \end{align}
\end{enumerate}
\end{thm}

The proof of the theorem can be found in
Section~\ref{SecProofThmUpper}. For constant pairs $(r,\beta)$,
Theorem~\ref{ThmUpper} can be applied with a target FDR proportional
to $\numobs^{-\kappastar}$ to show that both BH and BC achieve the
optimal decay of the combined FDR-FNR up to constant factors, as
stated formally below.

\begin{cor}\label{CorCoolUpper}
For $\beta < r$ and $\qstar = \cstar n^{-\kappastar}$ with $0 < \cstar
\leq \min\big\lbrace \frac{1}{24},~\frac{1}{6\Zlower}\big\rbrace$, the
BH and BC procedures with target FDR $\qstar$ satisfy
\begin{align}
\Risk_{n} \lesssim n^{-\kappastar} .
\end{align}
\end{cor}

\begin{figure}[t]
\centering
\begin{subfigure}{0.5\linewidth}
\includegraphics[width=\linewidth]{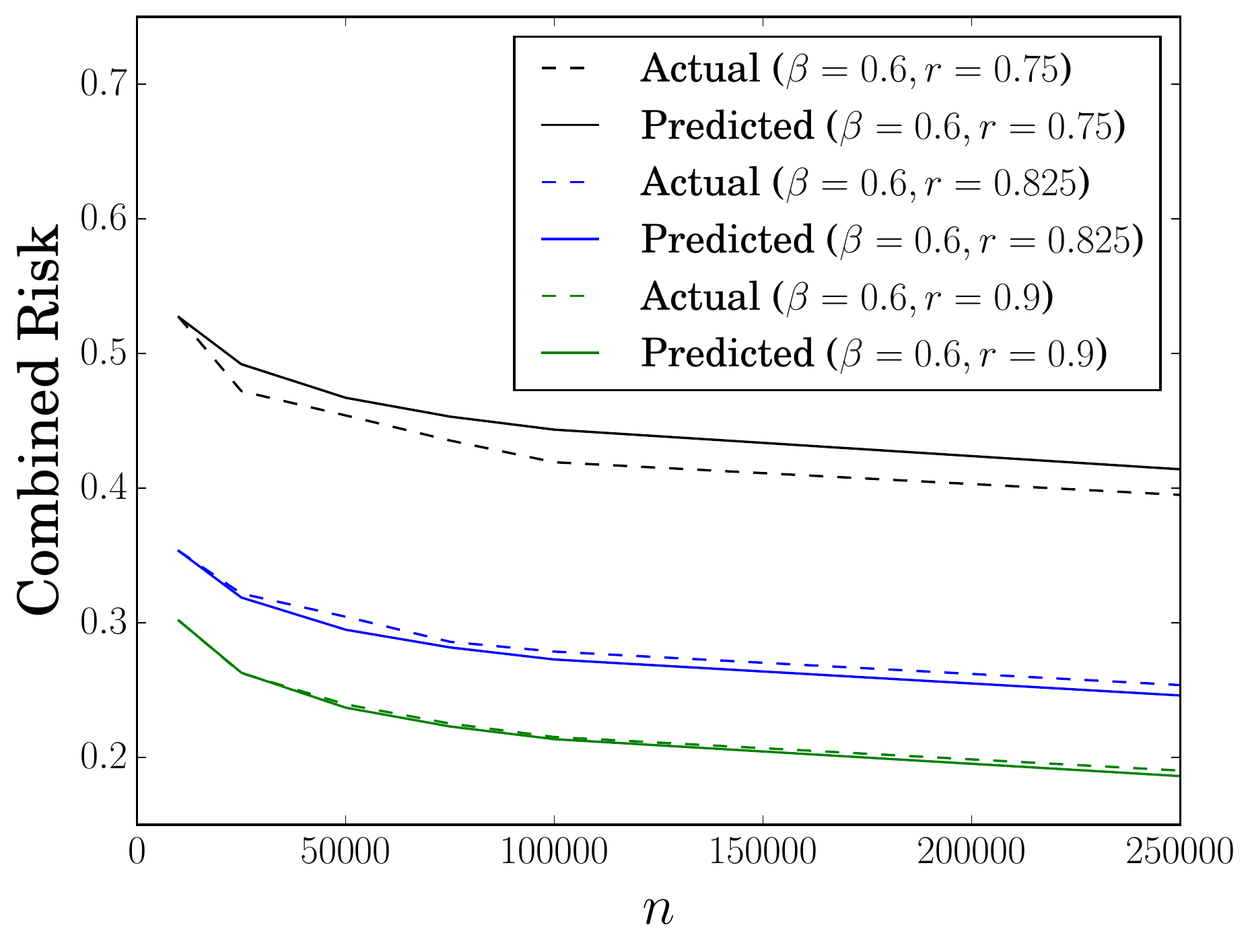}
\subcaption{$\gamma = 1$}
\end{subfigure}%
\begin{subfigure}{0.5\linewidth}
\includegraphics[width=\linewidth]{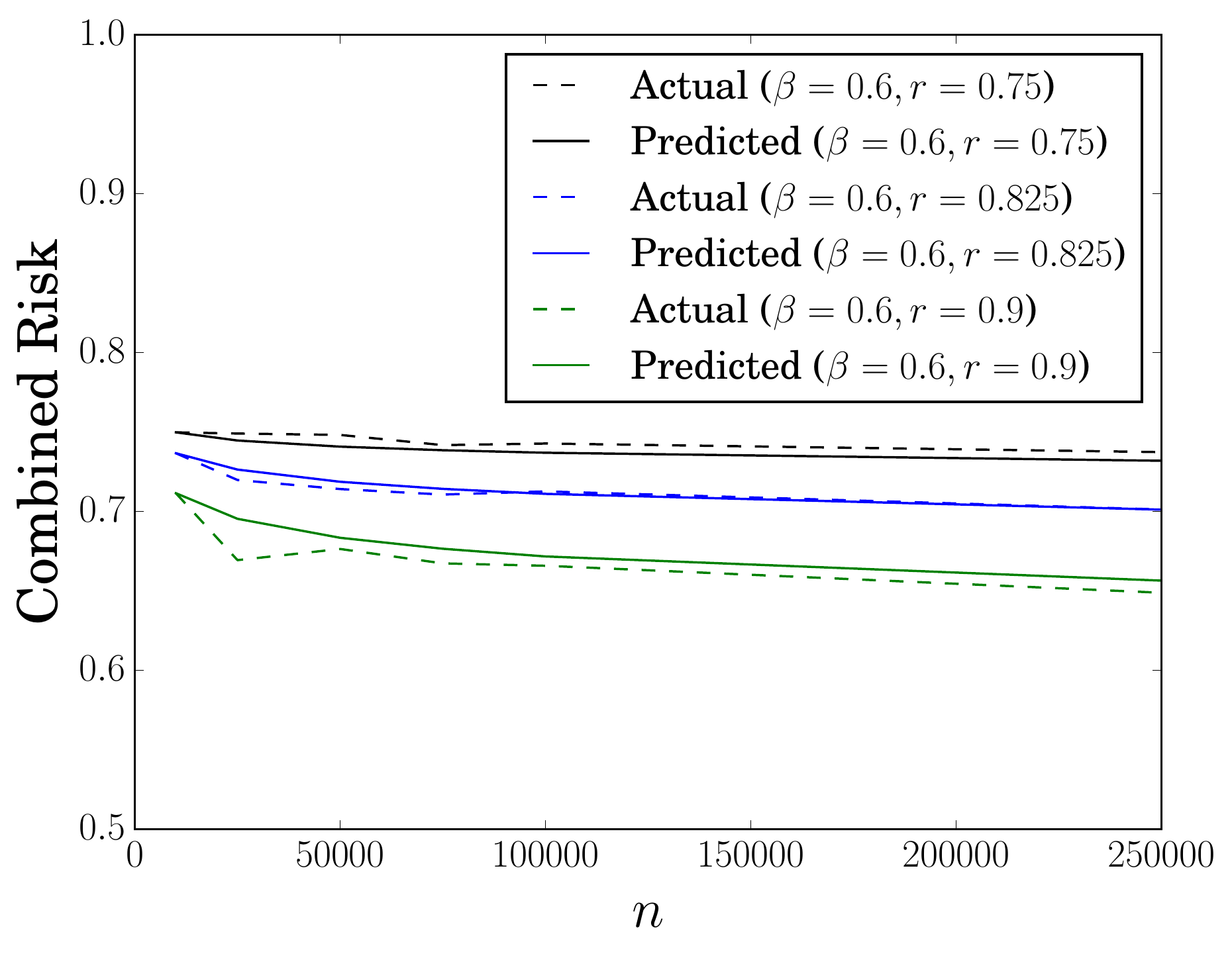}
\subcaption{$\gamma = 2$}
\end{subfigure}
\caption{Results of simulations comparing the predicted combined risk
  with the actual experimentally-observed risk for the BH
  procedure. Agreement is good across the board and improves as the
  gap $(r - \beta)$ increases. We believe the latter phenomenon arises
  because the sampling error is a smaller fraction of the risk as the
  separation increases.}
\label{FigUpper}
\end{figure}

To help visualize the result of Corollary~\ref{CorCoolUpper},
Figure~\ref{FigUpper} displays the results of some simulations of the
BH procedure that show correspondence between its performance and the
theoretically predicted rate of $\numobs^{-\kappanstar}$. Despite
the optimality, Figures~\ref{FigKappaStar} and~\ref{FigUpper} paint
a fairly dark picture from a practical point of view: while asymptotic
consistency can be achieved when $r > \beta$, the convergence of the
risk to zero can be extremely slow, exhibiting nonparametric rates
far slower than $n^{-1/2}$.  Figure~\ref{FigUpper} shows in particular
that the decay to zero may be barely evident even for sample sizes as
large as $n = 250,000$, even with comparatively strong signals.

The nonparametric nature may arise because the dimensionality 
of the decision space increases linearly with sample size, and
asymptotically, the upside of having increasing data seems to \emph{just}
overcome the downside of having to make an increasing number of decisions.
However, non-asymptotically, one cannot hope to drive both FDR and FNR to zero
at any practical sample size in this general setting, at least when the
mean signal lies below the maximum of the nulls (i.e., $r_{n} < 1$). 

\paragraph{Regime of linear sparsity:}

We turn to the regime of \emph{linear sparsity}---that is, when the number of
signals scales as $\pi_1 \numobs$ for some scalar $\pi_1 \in (0,1)$.  Recalling
that we have parameterized the number of signals as $\numobs^{1 - \betan}$, some
algebra leads to $\beta_\numobs = \frac{\log \frac{1}{\pi_1}}{\log
  n}$, so both Theorem~\ref{ThmLower} and Theorem~\ref{ThmUpper}
predict an upper and lower bound on the risk of the form
\begin{align}\label{EqLowerLinear}
c_{0} \exp\Big ( c_{1}\left[\frac{\log n}{\log
    \frac{1}{\pi_1}}\right]^{\frac{\gamma - 1}{\gamma}} \Big ) \cdot
n^{-\kappanstar}.
\end{align}
Note that here we overload the exponent $\kappanstar$ to the case when
it is nonconstant.  In order to interpret this result, observe that if
$r_{n} = r$ is constant, then $\kappanstar = \frac{r}{2^{\gamma}} -
o(1)$, so the rate is $n^{-r/2^{\gamma}}$ up to subpolynomial factors
in $n$.  On the other hand, if $r_{n} = \frac1{\log^{\frac{\gamma -
      1/2}{\gamma}}{n}}$ is at the extreme lower limit permitted by
the lower bound (ii) in~\eqref{EqnAnnoying2}, then it is not hard to
see that $\kappanstar \approx \log^{-\frac{\gamma - 1/2}{\gamma}}{n}$,
which ensures that $n^{\kappanstar} \gg \exp\big(\log^{\frac{\gamma -
    1}{\gamma}}{n}\big)$, so that the risk~\eqref{EqLowerLinear} still
approaches zero asymptotically, albeit subpolynomially in $n$.





\section{Proofs}
\label{SecProofs}

We now turn to the proofs of our main results, namely
Theorems~\ref{ThmLower} and~\ref{ThmUpper}, along with
Corollary~\ref{CorCool}.


\subsection{Proof of Theorem~\ref{ThmLower}}
\label{SecProofThmLower}
  
The main idea of the proof is to reduce the problem of lower bounding
the $\FNR_{n}$ of threshold-based procedures that use
random, data-dependent thresholds $T_n$, to the easier
problem of lower bounding the $\FNR_{n}$ of threshold-based procedures
that use a deterministic, data-independent threshold $t_{n}$.  We refer to
the latter class of procedures as \emph{fixed threshold procedures}, and
we parameterize them by their target FDR $q_\numobs =
\numobs^{-\kappan}$. Concretely, we define the \emph{critical threshold}, derived from
the critical regime boundary $\rmin$ from
equation~\eqref{eq:rn-minus}, by
\begin{align}
 \label{eq:tau-minus}
\taumin( \kappan ) & \defn \big ( \gamma \rmin\big ( \kappan \big )
\log{n} \big ) ^{1/\gamma} \quad \equiv \quad \taumin (q_\numobs )
\defn \left ( \gamma \rmin\left ( \frac{\log(1/q_\numobs)}{\log n}
\right ) \log{n} \right ) ^{1/\gamma}.
\end{align}
Here and throughout the proof, we express $\taumin$ and
$\rmin$ as functions of $q_\numobs$ rather than $\kappa_n$; this
formulation turns out to make certain calculations in the proof simpler
to express. 

\paragraph{From data-dependent threshold to fixed threshold:}

Our first step is to reduce the analysis from data-dependent to fixed
threshold procedures.  In particular, consider a threshold procedure,
using a possibly random threshold $T_\numobs$, that satisfies the FDR
uppper bound $\FDR_\numobs(T_\numobs) \leq q_\numobs$.  We claim that the
FNR of any such procedure must be lower bounded as
\begin{align}
  \label{eq:FNR-lower}
\E[\FNP_\numobs \big ( T_\numobs \big )] ~\geq~ \frac{\FNR_\numobs \big (
  \taumin \big ( 4q_\numobs \big ) \big ) }{16}.
\end{align}
This lower bound is crucial, as it reduces the study of random threshold procedures (LHS) to
study of fixed threshold procedures (RHS).
In order to establish the claim~\eqref{eq:FNR-lower}, define the
events
\begin{align*}
\event_{1} \defn \eventwrap{T_\numobs \geq \taumin\big (
  4q_\numobs\big )}, \quad \mbox{and} \quad \event_{2} \defn
\eventwrap{\FNP_\numobs\big ( \taumin\big ( 4q_\numobs \big )
  \big ) \geq \frac{\FNR_\numobs\big ( \taumin\big (
    4q_\numobs \big ) \big ) }{2}}.
\end{align*}
The following lemma guarantees that both of these events have
a non-vanishing probability:
\begin{lem}
  \label{LemEventBounds}
For any threshold $T_\numobs$ such that
    $\FDR_\numobs(T_\numobs) \leq q_\numobs$, we have
    \begin{subequations}
\begin{align}
  \label{EqnEventBounds}
  \P[\event_1] \; \stackrel{(a)}{\geq} \; 3/8, \quad \mbox{and} \quad
  \P[\event_2] \; \stackrel{(b)}{\geq} \; 3/4.
\end{align}
  \end{subequations}
\end{lem}
\noindent The proof of this lemma can be found in Appendix~\ref{AppLemEventBounds}.
Using this result, we now complete the proof of
claim~\eqref{eq:FNR-lower}.  Define the event
\begin{align*}
 \event & \defn \eventwrap{\FNP_\numobs(T_\numobs) \geq
   \frac{\FNR_\numobs \big ( \taumin(4 q_\numobs) \big )
   }{2}}.
\end{align*}
The monotonicity of the function $t \mapsto \FNP_{n}(t)$ ensures that
the inclusion $\event \supseteq \event_{1} \cap \event_{2}$ must hold.
Consequently, we have
\begin{align*}
\P[\event] & \geq \P[\event_1 \cap \event_2] \; \geq \; \P[\event_2] -
\P[\event_1^c] \; \stackrel{(i)}{\geq} \; \frac{3}{4} - \frac{5}{8} \;
= \; 1/8,
\end{align*}
where step (i) follows by applying the probability bounds from
Lemma~\ref{LemEventBounds}.

Finally, by Markov's inequality, we have
\begin{align*}
\FNR_\numobs(T_\numobs) \; = \; \Exs[ \FNP_\numobs(T_\numobs)] & \geq
\P[\event] \; \frac{\FNR_\numobs \big ( \taumin(4 q_\numobs)
  \big ) }{2} \; \geq \; \frac{\FNR_\numobs \big ( \taumin(4
  q_\numobs) \big ) }{16},
\end{align*}
which establishes the claim~\eqref{eq:FNR-lower}.

Our next step is to lower bound the FNR for choices of the threshold
$t \geq \taumin( q_\numobs)$:

\begin{lem}
\label{LemFNRFixed}
For any $t \geq \taumin(q_\numobs)$, we have
\begin{align}
  \label{EqnFNRFixed}
\FNR_\numobs (t) & \geq \begin{cases}
  \frac{\zeta^{2\beta_\numobs^{\frac{1 -
          \gamma}{\gamma}}}}{\Zlower} \cdot
  n^{-\dgamma{\beta_\numobs + \kappan}{r}} &~\text{if}~r >
  \rmin\big ( \kappan(q_\numobs) \big ) , \\ \frac{1}{2}
  &~\text{otherwise,}
\end{cases}
\end{align}
where $\zeta$ was previously defined~\eqref{EqZetaDef}.
\end{lem}
\noindent The proof of this lemma can be found in
Appendix~\ref{SecProofLemFNRFixed}. 
Armed with Lemma~\ref{LemFNRFixed} and the lower
bound~\eqref{eq:FNR-lower}, we can now complete the proof of
Theorem~\ref{ThmLower}.  We split the argument into two cases:

\paragraph{Case 1:}  First, suppose that
$r \leq \rmin(\kappa_n(4q_{\numobs}))$.  In this case, we have
\begin{align*}
\FNR_\numobs(T_\numobs) & \stackrel{(i)}{\geq} \frac{\FNR_\numobs\big
  ( \taumin\big ( 4q_\numobs \big ) \big ) }{16} \;
\stackrel{(ii)}{\geq} \frac{1}{32},
\end{align*}
where step (i) follows from the lower bound~\eqref{eq:FNR-lower}, and
step (ii) follows by lower bounding the FNR by $1/2$, as is guaranteed
by Lemma~\ref{LemFNRFixed} in the regime $r \leq
\rmin(\kappan(4q_{\numobs}))$.

\paragraph{Case 2:}  Otherwise, we may assume that
$r > \rmin( 4 q_\numobs)$.  In this case, we have
\begin{align*}
  \FNR_\numobs(T_\numobs) & \stackrel{(i)}{\geq}
  \frac{\FNR_\numobs\big ( \taumin\big ( 4q_n
    \big ) \big ) }{16} 
\; \stackrel{(ii)}{\geq} \frac{\big ( 4\zeta \big )
    ^{2\beta_\numobs^{\frac{1 - \gamma}{\gamma}}}}{\Zlower} \cdot
  n^{-\dgamma{\beta_\numobs + \kappan}{r}}.
\end{align*}
Here step (i) follows from the lower bound~\eqref{eq:FNR-lower},
whereas step (ii) follows from applying Lemma~\ref{LemFNRFixed} in the
regime $r > \rmin(\kappan(4 q_\numobs))$.
With some further algebra, we find that
\begin{align*}
\FNR_\numobs(T_\numobs) & \geq \frac{1}{\Zlower} \; \exp\big (
2\log\big ( 4\zeta \big ) \cdot \beta_\numobs^{\frac{1 -
    \gamma}{\gamma}} \big ) \; n^{-\dgamma{\beta + \kappa_n}{r}} \; =
\; c_0 \exp\big ( c_1 \beta_\numobs^{\frac{1 - \gamma}{\gamma}} \big )
\; n^{-\dgamma{\beta + \kappa_n}{r}},
\end{align*}
where $c_0 \defn \frac{1}{\Zlower}$ and $c_1 \defn 2\log\big (
4\zeta \big )$.  Note that since $\Zlower > 0$ and $\zeta \geq 1$,
both of the constants $c_0$ and $c_1$ are positive, as
claimed in the theorem statement.



\subsection{Proof of Corollary~\ref{CorCool}}
\label{SecProofCorCool}

We now turn the proof of Corollary~\ref{CorCool}.  Although it can be
proved from the statement of Theorem~\ref{ThmLower}, we instead prove
it more directly, as this allows us to reuse parts of the proof of
Lemma~\ref{LemFNRFixed}, thereby saving some additional messy
calculations.

First, we verify that there is indeed a unique solution
$\kappastar$ to the fixed point equation~\eqref{EqnCool}.  Define the
function as $g(\kappa) \defn \dgamma{\beta + \kappa}{r}^{1/\gamma} - \kappa^{1/\gamma}$.
Clearly the solutions to~\eqref{EqnCool} are the roots of $g$. We would like to argue that any such root must occur in $[0,~r - \beta)$ and that in fact $g$ has a unique root in this interval. For the first claim, note that $g(r - \beta) = -(r - \beta)^{1/\gamma} < 0$. On the other hand,
we have
\begin{align*}
g'(\kappa) = \begin{cases}
-\frac{1}{\gamma}\left[\left(\beta + \kappa\right)^{-\frac{\gamma - 1}{\gamma}} + \kappa^{-\frac{\gamma - 1}{\gamma}}\right] &~\text{if}~ 0 \leq \kappa < r - \beta, \\
\frac{1}{\gamma}\left[\left(\beta  + \kappa\right)^{-\frac{\gamma - 1}{\gamma}} - \kappa^{-\frac{\gamma - 1}{\gamma}}\right] &~\text{if}~ \kappa > r - \beta .  
\end{cases}
\end{align*}
It is immediately clear that $g'(\kappa) < 0$ for $0 \leq \kappa < r - \beta$ and, since $\beta + \kappa > \kappa$, we may also deduce that $g'(\kappa) < 0$ for $\kappa > r - \beta$, so $g$ is decreasing on its domain. Therefore, $g(\kappa) < g(r - \beta) < 0$ for all $\kappa > r - \beta$.\footnote{Note that we have suppressed the issue of non-differentiability of $g$ at $\kappa = r - \beta$. We may do so because it is left- and right-differentiable at this point, and we argue separately for the intervals $[0,~r - \beta)$ and $[r - \beta,~\infty)$.} We conclude that any root of $g$ must occur on $[0,~r - \beta)$. To finish the argument, note that $g(0) > 0 > g(r - \beta)$, so that $g$ does indeed have a root on $[0,~r - \beta)$. 

Turning now to the proof of the lower bound~\eqref{EqnCleanLower}, let
$\newdisc$ be an arbitrary threshold-based multiple testing procedure.
We may assume without loss of generality that
\begin{align}
 \FDR_\numobs(\newdisc) \leq \min \big \{ \numobs^{-\kappanstar},
 \frac{1}{24} \big \} \leq \cfun(\beta, \gamma) n^{-\kappanstar},
\end{align}
where the quantity $\cfun(\beta, \gamma) \geq 1$ was defined in the
statement of Theorem~\ref{ThmLower} (otherwise, the claimed lower
bound~\eqref{EqnCleanLower} follows immediately).

Applying the second part of Lemma~\ref{LemFNRFixed} and defining $\tilde{c} = 4\cfun(\beta, \gamma)$ , we
conclude that
\begin{align*}
\FNR_\numobs\big ( T_\numobs \big ) & \geq \frac{\FNR_\numobs \big (
  \taumin\big ( \tilde{c} n^{-\kappanstar} \big ) \big ) }{16} \\
& \geq \frac{\big (\tilde{c}\zeta \big ) ^{2\beta^{\frac{1 -
        \gamma}{\gamma}}}}{\Zlower} \cdot n^{-\dgamma{\beta +
    \kappa^{\ast}}{r}} \\
& = \frac{\big (\tilde{c}\zeta \big )
  ^{2\beta^{\frac{1 - \gamma}{\gamma}}}}{\Zlower} \cdot
n^{-\kappanstar} \\
& = c'n^{-\kappanstar} .
\end{align*}


\subsection{Proof of Theorem~\ref{ThmUpper}}\label{SecProofThmUpper}

We now aim to show that the Benjamini-Hochberg (BH) and
Barber-Cand\`{e}s (BC) algorithms achieve the minimax
rate~\eqref{EqnLower} when $r_{\numobs} >
\rmin(\kappan(c_{\alg}q_{\numobs}))$, where $\alg \in \lbrace
\BH,~\BC\rbrace$ and $c_{\alg}$ is the algorithm-dependent constant
defined in~\eqref{EqZetaDef}.  The proof strategy for both algorithms
is essentially the same. Given a target FDR rate $q_{n}$, we apply
each algorithm $q_{n}$ as the target FDR level and prove that the
resulting threshold satisfies $t_{\alg} \leq
\taumin(c_{\alg}q_{\numobs})$ with high probability.  The known
properties of the algorithms guarantee the required FDR bounds
~\cite[as studied by the authors of
][]{Ari16DistFreeFDR,Bar15Knockoffs,Ben95FDR}, while the following
converse to Lemma~\ref{LemFNRFixed} provides the requisite upper
bounds on the FNR.

\begin{lem}
\label{lem:FNR-fixed-converse}
If $r_{\numobs} > \rmin\big ( c q_n \big ) $ and $t
\leq \taumin(c q_n)$ for some $c > 0$, then we have
\begin{align*}
\FNR_\numobs\big ( t \big )  & \leq \frac{\big(\max \big\lbrace c,~1/c\big\rbrace \cdot \zeta\big)^{2\beta_\numobs^{\frac{1 -
        \gamma}{\gamma}}}}{\Zupper} \cdot n^{-\dgamma{\beta_\numobs +
    \kappan}{r}} ,
\end{align*}
where constant $\zeta$ is defined in~\eqref{EqZetaDef}. 
\end{lem}


\subsubsection{Achievability result for the BH procedure}

In this section, we prove that BH achieves the lower bound whenever
$r_\numobs > \rmin\big (c_{\BH}q_\numobs \big )$.  Specifically, we
prove the claim~\eqref{EqnBHBound} stated in Theorem~\ref{ThmUpper}.

As described in the previous section, the key step is  to compare $\tBH$, the random threshold set by BH with target FDR of $q_\numobs$, to the critical value
$\tauminBH := \taumin(c_{\BH} q_\numobs)$. To this end, we argue that
\begin{align}
  \label{EqBHThreshBd}
\P\big(\tBH > \tauminBH \big ) \leq \exp\big ( -\frac{n^{1
    -\rmax}}{24} \big ).
\end{align}

We begin by showing how to derive the upper bound~\eqref{EqnBHBound}
from the probability bound~\eqref{EqBHThreshBd}.  Note that since BH
is a valid FDR control procedure, we necessarily have $\FDR_{n}(\tBH)
\leq q_\numobs$. To bound the FNR, first let $\event = \eventwrap{\tBH
  \leq \tauminBH}$ and let $\FNR_{n}(\cdot~|~\event)$ and
$\FNR_{n}(\cdot~|~\event^{c})$ denote the $\FNR_{n}$ conditional on
the event and its complement, respectively. In this notation, the
bound~\eqref{EqBHThreshBd}, together with
Lemma~\ref{lem:FNR-fixed-converse}, implies that
\begin{align*}
\FNR_{n}(\tBH) & \leq \P\big(\event\big) \cdot
\FNR_{n}\big(\tauminBH~|~\event\big) + \P\big(\event^{c}\big) \\ &
\leq \FNR_{n}\big(\tauminBH\big) + \P\big(\event^{c}\big) \\ & \leq
\FNR_{n}\big(\tauminBH\big) + \exp\big ( -\frac{n^{1 - \rmax}}{24}
\big ) \\ & \leq \frac{\zeta_{\BH}^{2\beta_\numobs^{\frac{1 -
        \gamma}{\gamma}}}}{\Zupper} \cdot n^{-\dgamma{\beta_\numobs +
    \kappan}{r_\numobs}} + \exp\big ( -\frac{n^{1 - \rmax}}{24} \big )
\\ & \leq \frac{2\zeta_{\BH}^{2\beta_\numobs^{\frac{1 -
        \gamma}{\gamma}}}}{\Zupper} \cdot n^{-\dgamma{\beta_\numobs +
    \kappan}{r_\numobs}},
\end{align*}
where the final step uses the definition~\eqref{EqNminUpper} of
$\nminu$, and the fact that $\frac{1}{\Zupper n} \leq
\frac{\zeta_{\BH}^{2\beta_\numobs^{\frac{1 -
        \gamma}{\gamma}}}}{\Zupper} \cdot n^{-\dgamma{\beta_\numobs +
    \kappan}{r_\numobs}}$, which is easily verified by noting that
$\zeta_{\BH}^{2\beta_\numobs^{\frac{1 - \gamma}{\gamma}}} \geq 1$ and
$\dgamma{\beta_\numobs + \kappan}{r_\numobs} \leq 1$.

We now prove the bound~\eqref{EqBHThreshBd} with an argument using
$p$-values and survival functions that parallels that
of~\citet{Ari16DistFreeFDR} but sidesteps CDF asymptotics.  We study
the relationship between the population survival function $\survival$
and the empirical survival function $\empsurv$, defined by
\begin{align}\label{EqEmpSurvDef}
\empsurv\big ( t \big ) & = \big ( 1 - \frac{1}{n^{\beta_\numobs}}
\big ) \cdot \empsurvnull\big ( t \big ) + \frac{1}{n^{\beta_\numobs}}
\cdot \empsurvsig \big ( t \big ),\\ \text{where ~}~ \quad
\empsurvnull\big ( t \big ) &= \frac{1}{n - n^{1 -
    \beta_\numobs}}\sum_{i \in \nulls} \one\big ( X_{i} \geq t \big )
~~\text{and} ~~ \empsurvsig\big ( t \big ) = \frac{1}{n^{1 -
    \beta_\numobs}}\sum_{i \notin \nulls} \one\big ( X_{i} \geq t \big
) .\nonumber
\end{align}
Now, sort the observations in \emph{decreasing} order, 
so that $X_{(1)} \geq X_{(2)} \geq \cdots \geq X_{(\numobs)}$, and
define $p$-values
\begin{align}
  \label{EqBHPValues}
p_{(i)} = \survival\left(X_{(i)}\right) ~~ \text{and} ~~
\empsurv\left(X_{(i)}\right) = \frac{i}{n},
\end{align}
so that $p_{(1)} \leq p_{(2)} \leq \cdots \leq p_{(\numobs)}$ are in
\emph{increasing} order.  Then, we may characterize the indices
rejected by BH as those satisfying $X_{i} \geq X_{(\iBH)}$, where
  \begin{align}
    \label{EqIBHDef}
\iBH = \max \big \lbrace 1 \leq i \leq n \colon
\survival\left(X_{(i)}\right) \leq
q_\numobs\empsurv\left(X_{(i)}\right) \big \rbrace.
\end{align}
Moving $\tBH$ within $\big(X_{(\iBH + 1)},~X_{(\iBH)}\big]$ if
      necessary, we may therefore assume $\survival(\tBH) =
      q_{\numobs}\empsurv(t)$ whenever $t < \tBH$, and combining this
      knowledge with~\eqref{EqEmpSurvDef}, we obtain the chain of
      inclusions
\begin{align}
\event^{c} = \eventwrap{\tBH \leq \tauminBH}& \supset
\eventwrap{\survival\big ( \tauminBH \big ) \leq q_\numobs\empsurv\big
  ( \tauminBH \big )} \nonumber \\ & \supset \eventwrap{\survival\big
  ( \tauminBH \big ) \leq \frac{q_\numobs}{n^{\beta_\numobs}} \cdot
  \frac{W_\numobs}{n^{1 - \beta_\numobs}}} = \colon \widetilde{\event}
,
\end{align}
where $W_\numobs = \sum_{i \notin \nulls} \one\big(X_{i} \geq
\tauminBH\big) \sim \mathrm{Bin}\big ( \survival\big ( \tauminBH -
\mu_{\numobs} \big ) ,~n^{1 - \beta_\numobs} \big ) $.

We now argue that $\survival\big ( \tauminBH \big ) \leq
\frac{q_\numobs}{4 \numobs^{\beta_\numobs}}$, so that
$\P\left(\event\right) \geq \P\left(W_\numobs >\frac{n^{1 -
    \beta_\numobs}}{4}\right)$.  For this, observe that by the
definition of $\rmin$ in~\eqref{eq:rn-minus} and the upper tail
bound~\eqref{EqTailBounds}, we have
\begin{align*}
\log\survival\big ( \tauminBH \big ) & \leq
-\rmin\left(c_{\BH}q_{n}\right) \log{n} + \log \frac{1}{\Zupper} \\ &
\leq -\beta_\numobs \log{n} + \log (c_{\BH}q_\numobs) -
\log\frac{1}{6\Zlower} + \log \frac{1}{\Zupper} \\ & = -\beta_\numobs
\log{n} + \log{q_{n}} + \log \frac{6c_{\BH}\Zlower}{\Zupper} \\ & =
\log{\frac{q_\numobs}{6n^{\beta_\numobs}}} ~<~
\log{\frac{q_\numobs}{4n^{\beta_\numobs}}}.
\end{align*}
We conclude 
\begin{align*}
\P\left(\tBH > \tauminBH \right) ~\leq~ 1 -
\P\left(\widetilde{\event}\right) ~\leq~ 1 - \P\left(W_\numobs
>\frac{n^{1 - \beta_\numobs}}{4}\right) ~=~ \P\left(W_\numobs \leq
\frac{n^{1 - \beta_\numobs}}{4}\right).
\end{align*}
Finally, by a Bernstein bound, we find
\begin{align*}
\P\left(W_\numobs \leq \frac{n^{1 - \beta_\numobs}}{4}\right) & \leq
\P\left(W_\numobs \leq \frac{\E\left[W_\numobs\right]}{2}\right) \\ &
\leq \exp\big(-\frac{\E\left[W_\numobs\right]}{12}\big) \\ & \leq
\exp\big(-\frac{n^{1 - \beta_\numobs}}{24}\big) \\ & \leq
\exp\big(-\frac{n^{1 - \rmax}}{24}\big),
\end{align*}
where we have used the fact that $\tauminBH \leq \mu_{\numobs}$ to
conclude that $\survival\left(\tauminBH -\mu_{\numobs}\right) \geq
\frac{1}{2}$ and therefore $\E[W_{n}] \geq \frac{n^{1 -
    \beta_n}}{2}$. We have therefore established the required
claim~\eqref{EqBHThreshBd}, concluding the proof of optimality of the
BH procedure.


\subsubsection{Achievability result for the BC procedure}
\label{subsec:bc}

Our overall strategy for analyzing BC procedure resembles the one we
used for the BH procedure. As with our analysis of the BH procedure,
we define $\tauminBC := \taumin(c_{\BC}q_{n})$ and derive the
bound~\eqref{EqnBHBound} by controlling the algorithm's threshold as
\begin{align}
  \label{eq:BC-thresh-bd}
\P\left(\tBC > \tauminBC \right) \leq q_\numobs + \exp\big (
-\frac{n^{1 - \rmax}}{24} \big ).
\end{align}
Since the proof of equation~\eqref{EqnBCBound} from the
bound~\eqref{eq:BC-thresh-bd} is essentially identical to the
corresponding derivation for the BH procedure, we omit it.
\noindent We now prove the bound~\eqref{eq:BC-thresh-bd} by an
argument somewhat different than that used in analyzing the BH
procedure. Define the integers
\begin{align*}
N_{+} \big ( t \big ) = \sum_{i = 1}^{n} \one\big ( X_{i} \geq t \big
) ~~ \text{and} ~~ N_{-}\big ( t \big ) = \sum_{i = 1}^{n} \one\big (
X_{i} \leq -t \big ) .
\end{align*}
Then, the definition of the BC procedure gives
\begin{align*}
\tBC = \inf \big\lbrace t \in \R \colon \frac{1 + N_{-}(t)}{1 \vee
  N_{+}(t)} \leq q_\numobs\big\rbrace .
\end{align*}
To prove~\eqref{eq:BC-thresh-bd}, it therefore suffices to show that
\begin{align}
  \label{eq:bc-ratio-tau-upper}
\P\Bigg(\frac{1 + N_{-}(\tauminBC)}{1 \vee
  N_{+}(\tauminBC)} > q_\numobs\Bigg) \leq q_\numobs + \exp\big ( -\frac{n^{1 - \rmax}}{24}
\big ) .
\end{align}
We prove the bound~\eqref{eq:bc-ratio-tau-upper} in two parts:
\begin{subequations}
\begin{align}
\label{eq:bc-ratio-tau-upper-plus}
\P\big(1 \vee N_{+}(\tauminBC) < \frac{n^{1 -
    \beta_\numobs}}{4}\big) \leq \exp\big ( -\frac{n^{1 - \rmax}}{24}
\big ),\\
\label{eq:bc-ratio-tau-upper-minus} 
\P\Bigg(1 + N_{-}(\tauminBC) > q_{\numobs} \cdot
\frac{n^{1 - \beta_{n}}}{4}\Bigg) \leq q_\numobs.
\end{align}
\end{subequations}
These bounds are a straightforward consequence of elementary Bernstein
bounds, and together they imply the
claim~\eqref{eq:bc-ratio-tau-upper}. We explain them below.

The lower bound~\eqref{eq:bc-ratio-tau-upper-plus} follows because $1
\vee N_{+}(\tauminBC) \geq N_{+}(\tauminBC)$ and
$N_{+}(\tauminBC)$ is the sum of two binomial random
variables, corresponding to nulls and signals, respectively, and the
latter has a $\survival\big ( \tauminBC - \mu \big ) \geq
\frac{1}{2}$ probability of success. More precisely, we may write
$N_{+}\left(\tauminBC\right) = N_{+}^{\nulltxt} +
N_{+}^{\sigtxt}$, with
\begin{align*}
N_{+}^{\nulltxt} \sim
\Bin\left(\survival\left(\tauminBC\right),~n - n^{1 -
  \beta_{n}}\right) ~~ \text{and} ~~ N_{+}^{\sigtxt} \sim
\Bin\left(\survival\left(\tauminBC - \mu_{n}\right),~n^{1 -
  \beta_{n}}\right) ,
\end{align*}
implying $N_{+}\left(\tauminBC\right) \geq N_{+}^{\sigtxt}$,
whence
\begin{align*}
\E\big[N_{+}(\tauminBC)\big] \geq \E\big[N_{+}^{\sigtxt}\big]
= n^{1 - \beta_\numobs} \cdot \survival\big ( \tauminBC - \mu_{n}
\big ) \geq \frac{n^{1 - \beta_n}}{2},
\end{align*}
where we have used the fact that $\tauminBC \leq \mu_{n}$. With
this bound in hand, a Bernstein bound yields
\begin{align*}
\P\big(N_{+}(\tauminBC) < \frac{n^{1 - \beta_n}}{4}\big) & \leq \P\big(N_{+}(\tauminBC) \leq \frac{\E\left[N_{+}(\tauminBC)\right]}{2}\big) \\
& \leq \exp\big ( -\frac{n^{1 - \beta_n}}{24} \big ) \leq \exp\big (
-\frac{n^{1 - \rmax}}{24} \big ) ,
\end{align*}
as required to prove equation~\eqref{eq:bc-ratio-tau-upper-plus}.
The proof of equation~\eqref{eq:bc-ratio-tau-upper-minus} follows a
similar pattern. Here, we note that $N_{-}(\tauminBC)$ is a
sum of two binomial random variables, with a total of $n$ trials, such
that---using the definition~\eqref{eq:rn-minus} of $\rmin$ and the upper bound on the tail~\eqref{EqTailBounds}---each one has probability of success upper bounded by $1 - \survival\big
( -\tauminBC \big ) \leq \frac{6\Zlower}{\Zupper} \cdot
c_{\BC}q_\numobs \numobs^{-\beta_n} = \frac{1}{8} \cdot q_{\numobs}\numobs^{-\beta_{\numobs}}$. Formally, we may write
$N_{-}(\tauminBC) = N_{-}^{\nulltxt} + N_{-}^{\sigtxt}$, with
\begin{align*}
N_{-}^{\nulltxt} \sim \Bin\left(1 -
\survival\left(-\tauminBC\right),~n - n^{1-
  \beta_{\numobs}}\right) ~~ \text{and} ~~ N_{-}^{\sigtxt} \sim
\Bin\left(1 - \survival\left(-\tauminBC - \mu\right),~n^{1 -
  \beta_{\numobs}}\right) .
\end{align*}
Since $1 - \survival\left(-\tauminBC - \mu\right) \leq 1 -
\survival\left(-\tauminBC\right)$, we deduce
\begin{align*}
\E\big[N_{-}(\tauminBC)\big] & \leq \big[1 -
  \survival\left(-\tauminBC\right) \big] \cdot n \leq 
\frac{q_{\numobs}}{2} \cdot \frac{n^{1 -
    \beta_n}}{4}.
\end{align*}
On the other hand, using the lower bound in~\eqref{EqTailBounds}, we find 
$1 - \survival\big( -\tauminBC \big ) \geq 6c_{\BC}q_{n}n^{-\beta_{n}}$. Using the additional fact that
$n - n^{1 - \beta_{n}} \geq \frac{n}{2}$ by~\eqref{EqnAnnoying1}, we may conclude that
\begin{align*}
\E\big[N_{-}(\tauminBC)\big] & \geq \E\big[N_{-}^{\nulltxt}\big] \\
                                           & = \big(n - n^{1 - \beta_{n}}\big) \cdot \big[1 - \survival\big( -\tauminBC \big )\big] \\
                                           & \geq \frac{n}{2} \cdot 6c_{\BC}q_{n}n^{-\beta_{n}} \\
                                           & \geq 3c_{\BC}q_{n} n^{1 - \beta_{n}} .
\end{align*}
By a Bernstein bound, it follows that
\begin{align*}
\P\Bigg(N_{-}(\tauminBC) \geq 
     q_\numobs \cdot \frac{n^{1 -
    \beta_n}}{4}\Bigg) & \leq \P\Bigg(N_{-}(\tauminBC) \geq 2\E\big[N_{-}(\tauminBC)\big] \Bigg) \\
    & \leq \exp\Bigg( -\frac{\E\big[N_{-}(\tauminBC)\big]}{4}\Bigg) \\
    & \leq \exp\big(-\frac{3c_{\BC}}{4} \cdot q_{n}n^{1 - \beta_{n}} \big) \\
    & \leq q_{n},
\end{align*}
where we have invoked the decay condition~\eqref{eq:qn-decay} for the last step. 


\subsection{Proof of Corollary~\ref{CorCoolUpper}}

The corollary is a nearly immediate consequence of
Theorem~\ref{ThmUpper}. We will prove it for both algorithms
simultaneously.  Observe that
\begin{align}
  \label{EqRMinCoolProof}
r_{\min}(\kappan(c_{\alg}\qstar)) = \beta + \kappastar + \frac{\log
  \frac{1}{6\cstar c_{\alg}\Zlower}}{\log{n}}.
\end{align}
Suppose for now that the decay condition~\eqref{eq:qn-decay} holds for
$\qstar$ and some choice of $n_{\min,\BC}$. Then,
using~\eqref{EqRMinCoolProof} and the fact that $r > \beta +
\kappastar$, we may choose $\nmin' \geq n_{\min,\BC}$ large enough so
that $r > r_{\min}(\kappan(c_{\alg}\qstar))$ for all $n \geq \nmin'$
and $\alg \in \left\lbrace \BH,~\BC\right\rbrace$. From
Theorem~\ref{ThmUpper}, we conclude that there exists a constant $c'$
such that both algorithms satisfy
\begin{align*}
n \geq \nmin'  \Longrightarrow \Risk_{n} \leq c'n^{-\kappastar} .
\end{align*}
By replacing $c'$ by $\tilde{c} = \max\left\lbrace
c',~\left(\nmin'\right)^{\kappastar}\right\rbrace$ (and recalling
$\Risk_{n} \leq 1$ always), we obtain $\Risk_{n} \leq
\tilde{c}n^{-\kappastar}$ for all $n \geq 1$, obtaining the claimed
result.

In order to check the decay condition~\eqref{eq:qn-decay}, note that,
as $\kappastar \leq r - \beta \leq 1 - \beta$, we have for
sufficiently large $n$ that
\begin{align*}
\frac{q_{\numobs}}{\log \frac{1}{q_{\numobs}}} =
\frac{n^{-\kappastar}}{\kappastar \log{n}} \geq \frac{4}{3c_{\BC}}
\cdot n^{-(1 - \beta)},
\end{align*}
which completes the proof.


\section{Discussion}
\label{SecDisc}

Despite considerable interest in multiple testing with false discovery
rate (FDR) control, there has been relatively little understanding of
the non-asymptotic trade-off between controlling FDR and the
analogous measure of power known as the false non-discovery rate
(FNR).  In this paper, we explored this issue in the context of the
sparse generalized Gaussians model, and derived the first
non-asymptotic lower bounds on the sum of FDR and FNR.  We
complemented these lower bounds by establishing the non-asymptotic
minimaxity of both the Benjamini-Hochberg (BH) and Barber-Cand\`{e}s
(BC) procedures for FDR control.  The theoretical predictions are
validated in simple simulations, and our results recover recent
asymptotic results~\cite{Ari16DistFreeFDR} as special cases.  Our work
introduces a simple proof strategy based on reduction to deterministic
and data-oblivious procedures.  We suspect this core idea may apply to
other multiple testing settings: in particular, since our arguments do
not depend on CDF asymptotics in the way that many classical analyses
of both global null testing and FDR control procedures do, we hope
they will be possible to adapt for other problems described below.

As mentioned after the statement of Theorem~\ref{ThmUpper}, the
practical implications of our results are somewhat pessimistic.  Even
for rather simple problems having $r - \beta$ of constant order, the
resulting rate at which the risk tends to zero can be far slower than
$n^{-1/2}$.  (Indeed, it seems like such a parametric rate is only
achievable when $\gamma = 1, r_n \to 1, \beta_n \to 0$.) Hence, in
practice, one must carefully consider whether good FDR or good FNR is
more important, as achieving both may not be possible unless most of
the signals to be identified are rather large.

\subsection*{Future work} A large part of the multiple testing
literature focuses on the development of valid FDR control procedures
that can gain power or precision by explicitly using prior knowledge
such as null-proportion adaptivity~\cite{Storey02,Storey04}, groups or
partitions of hypotheses~\cite{Bar15PFilter,hu2010false}, prior or
penalty weights~\cite{BH97,Gen06FDRWeights}, or other forms of
structure~\cite{Li16Struct,ramdas2017unified}, and it would be of
great interest to extend our techniques to such structured settings.
It is also important to handle the cases of positive or arbitrary
dependence~\cite{BY01,blanchard2008two,ramdas2017unified}, hence this
is another natural direction in which to extend our work.  (Such
extensions were explored
by~\cite{hall2008properties,Jin14RareWeak,hall2010innovated} for the
higher criticism statistic for the detection problem.)  Lastly, all of
the discussed literature is in the offline setting, and it could be of
interest to develop tight lower and upper bounds for online FDR
procedures~\cite{Fos08Alpha,Jav15Online}.

A general proof technique for establishing non-asymptotic lower bounds
in multiple testing remains an important direction for future work. As
our arguments are based on analytical calculations, they are sensitive
to the particular observation model under consideration. It would be
desirable to build on our approach to identify the key properties of
multiple testing problems that make them difficult or easy, and
establish lower bounds in terms of these properties. We hope our work
will help lay the foundation for progress in understanding the limits
of multiple testing more generally.


\subsection*{Acknowledgements}

This work was partially supported by Office of Naval Research Grant
DOD-ONR-N00014, Air Force Office of Scientific Research Grant
AFOSR-FA9550-14-1-0016, and Office of Naval Research grant W911NF-16-1-0368NSF.
In addition, MR was supported by an NSF Graduate Research Fellowship and a
Fannie and John Hertz Foundation Google Fellowship.


\appendix

\section{Proof of Lemma~\ref{LemEventBounds}}
\label{AppLemEventBounds}

This appendix is devoted to the proof of Lemma~\ref{LemEventBounds},
which was involved in the proof of Theorem~\ref{ThmLower}.  Our proof
makes use of the following auxiliary lemma:
\begin{lem}
\label{LemUniformThreshold}
For $q_{\numobs} \in (0,1/24)$, we have
\begin{align}
\label{eq:tau-FDP-lower}
\P \Big[ \FDP_\numobs\big ( t \big ) \geq 8 q_{n} \quad \mbox{for all $t
    \in [0, \taumin(4q_{n})]$} \Big] & \geq \frac{1}{2}.
\end{align}
\end{lem}
\noindent We return to prove this claim in
Appendix~\ref{AppLemUniformThreshold}.  For the moment, we take it as
given and complete the proof of Lemma~\ref{LemEventBounds}.


\paragraph{Control of $\event_1$:}  Let us now prove the first bound
in Lemma~\ref{LemEventBounds}, namely that $\P[\event_{1}] \geq
\frac{3}{8}$ where $\event_1 \defn \big \{ T_{\numobs} \geq
\taumin(4 q_\numobs) \big \}$.  So as to simplify notation,
let us define the event
\begin{subequations}
\begin{align}
  \label{EqnDefnEvent}
\Event & \defn \Big \{ \FDP_\numobs\big ( t \big ) \geq 8 q_\numobs
\quad \mbox{for all $t \in [0, \taumin(4q_\numobs)]$} \Big
\}.
\end{align}
Now observe that 
\begin{align}
  \label{EqnMailani}
\P \Big[ T_{\numobs} \geq \taumin \big( 4 q_\numobs \big)
  \Big] \; \geq \; \P \Big[ T_{\numobs} \geq \taumin\big (
  4q_\numobs \big ) \text{ and } \Event \Big] \; = \; \P [\Event] - \P
\big[ T_{\numobs} \leq \taumin\big ( 4q_\numobs \big ) \text{
    and } \Event \big].
\end{align}
\end{subequations}
Now by the definition~\eqref{EqnDefnEvent} of the event $\Event$, we
have the inclusion
\begin{align*}
\Big \{ T_{\numobs} \leq \taumin\big ( 4q_\numobs \big )
\text{ and } \Event \Big \} & \subseteq \Big \{
\FDP_{\numobs}(T_{\numobs}) \geq 8q_{\numobs} \Big \}.
\end{align*}
Combining with our earlier bound~\eqref{EqnMailani}, we see that
\begin{align*}
\P \Big[ T_{\numobs} \geq \taumin \big( 4 q_\numobs \big)
  \Big] & \geq \P [ \Event] - \P \big[ \FDP_{\numobs}(T_{\numobs})
  \geq 8q_{\numobs} \big].
\end{align*}
It remains to control the two probabilities on the right-hand side of
this bound.  Applying Lemma~\ref{LemUniformThreshold} guarantees that 
$\P [\Event] \geq \frac{1}{2}$.  On the other hand, by Markov's inequality, the assumed lower bound
$\FDR_{\numobs}\big(T_{\numobs}\big)\leq q_{n}$ implies that $\P \Big[
  \FDP_{\numobs}\big(T_{\numobs}\big) \geq 8 q_{\numobs} \Big] \leq
\frac{1}{8}$.  Putting together the pieces, we conclude that
\begin{align*}
\P[\event_1] \; = \; \P \Big[ T_{\numobs} \geq \taumin \big(
  4 q_\numobs \big) \Big] \; \geq \; \; \frac{1}{2} - \frac{1}{8} \; =
\; \frac{3}{8},
\end{align*}
as claimed.


\paragraph{Control of $\event_2$:}
Let us now prove the lower bound $\P [\event_2] \geq 3/4$.  We split
our analysis into two cases.

\paragraph{Case 1:} First, suppose that $r_\numobs > \rmin$. In
this case, we can write
\begin{align*}
\FNP_\numobs(t) = \frac{F_\numobs(t)}{n^{1 - \beta_\numobs}}, \qquad
\mbox{where $F_\numobs\big ( t \big ) \sim \mathrm{Bin}\big ( 1 -
  \survival\big ( t - \mu \big ) ,~n^{1 - \beta_\numobs} \big )$.}
\end{align*}
Since $\rmin > \beta_\numobs$, we have
\begin{align*}
\mu - \taumin \geq \big ( \gamma \log n \big )
^{1/\gamma}\big[r_\numobs^{1/\gamma} - \beta_\numobs^{1/\gamma}\big] \geq
\big ( \gamma \log{n} \big ) ^{1/\gamma} \cdot \big ( r_\numobs -
\beta_\numobs \big ) ^{1/\gamma},
\end{align*}
from which it follows that
\begin{align}
  \label{EqnSurviveBound}
\Exs[F_\numobs] \; = \; 1 - \survival\big ( t - \mu \big ) \geq
\frac{n^{\beta_\numobs - r_\numobs}}{\Zlower}.
\end{align}
Now by applying the Bernstein bound to the binomial random variable
$F_\numobs$, we have
\begin{align}
1 - \P[\event_{2}] = \P \left[ F_\numobs \leq \frac{\E [F_\numobs]}{2}
  \right] & \leq \exp\big ( -\frac{\E[F_\numobs]}{12} \big ) \nonumber
\\
& \stackrel{(i)}{\leq} \exp \big ( -\frac{n^{1 -
    r_\numobs}}{12\Zlower} \big ) \nonumber \\
\label{EqnCaseOne}
& \stackrel{(ii)}{\leq} \exp \big ( -\frac{n^{1 - \rmax}}{12\Zlower}
\big ),
\end{align}
where step (i) follows from the lower bound~\eqref{EqnSurviveBound},
and step (ii) follows since $r_\numobs < \rmax$ by assumption.

\paragraph{Case 2:}  Otherwise, we may assume that
$r_n \in \big( \beta_n, r_\numobs^{-1} \big)$.  In this regime, we
have the lower bound $\taumin - \mu \geq 0$, so that the
binomial random variable $F_\numobs$ stochastically dominates a second
binomial distributed as $\Ftil_\numobs \sim \mathrm{Bin}\big (
\frac{1}{2},~n^{1 - \beta_\numobs} \big ) $. By this stochastic
domination condition, it follows that
\begin{align*}
1 - \P \big[\event_{2}\big] \; \leq \; \P \Big[F_\numobs \leq
  \frac{n^{1 - \beta_\numobs}}{4} \Big ] & \leq \P \Big[\Ftil_\numobs
  \leq \frac{\E\left[\Ftil_\numobs\right]}{2} \Big].
\end{align*}
By applying the Bernstein bound to $\Ftil_\numobs$, we find that
\begin{align}
  \label{EqnCaseTwo}
  1 - \P \big[\event_{2}\big] & \leq \exp\big ( -\frac{n^{1 -
      \beta_\numobs}}{24} \big ) \; \leq \exp\big ( -\frac{n^{1 -
      \rmax}}{24} \big ),
\end{align}
where the final step follows since $\rmax > \beta_\numobs$. \\

\vspace*{.1in}
\noindent Putting together the two bounds~\eqref{EqnCaseOne}
and~\eqref{EqnCaseTwo}, we conclude that $\P [\event_2] \geq
\frac{3}{4}$ for all sample sizes $\numobs$ large enough to ensure
that
\begin{align}
\max\big\lbrace \exp\big (-\frac{n^{1 - \rmax}}{24\Zlower} \big ) ,~
\exp\big ( -\frac{n^{1 - \rmax}}{24} \big ) \big\rbrace \leq
\frac{1}{4} \label{EqNminLowerApplied},
\end{align}
as was claimed. Note that condition~\eqref{EqNminLowerApplied} is
identical to condition~\eqref{EqNminLower}, so that our definition of
$\nmin$ guarantees that~\eqref{EqNminLowerApplied} is satisfied.  This
completes the proof.

\subsection{Proof of Lemma~\ref{LemUniformThreshold}}
\label{AppLemUniformThreshold}

It remains to prove our auxiliary result stated in
Lemma~\ref{LemUniformThreshold}.  For notational economy, let $\tau =
\taumin(s)$ and let $\beta = \beta_\numobs$. The FDP at a
threshold $t$ can be expressed in terms of two binomial random
variables
\begin{align*}
L_\numobs\big ( t \big ) = \sum_{i \in \nulls} \one\big ( X_{i} \geq t
\big ) ~~ \text{and} ~~ W_\numobs\big ( t \big ) = \sum_{i \notin
  \nulls} \one\big ( X_{i} \geq t \big ) \leq n^{1 - \beta}.
\end{align*}
Here $L_\numobs(t)$ and $W_\numobs(t)$ correspond (respectively) to
the number of nulls, and the number of signals that exceed the
threshold $t$.  In terms of these two binomial random variables, we
have the expression
\begin{align*}
\FDP_\numobs(t) = \frac{L_\numobs(t)}{L_\numobs(t) + W_\numobs(t)}
\geq \frac{L_\numobs(t)}{L_\numobs(t) + n^{1 - \beta}} .
\end{align*}
Note that the inequality here follows by replacing $W_\numobs(t)$ by
the potentially very loose upper bound $n^{1 - \beta}$; doing so
allows us to reduce the problem of bounding the FDP to control of
$L_\numobs(t)$ uniformly for $t \in [0, \tau]$.
By definition of $L_\numobs(t)$, we have the lower bound
\begin{align*}
  \frac{L_\numobs(t)}{L_\numobs(t) + n^{1 - \beta}} \geq
  \frac{L_\numobs(\tau) }{L_\numobs(\tau) + n^{1 - \beta}} \qquad
  \mbox{for all $t \in [0, \tau]$.}
\end{align*}
Moreover, observe that
\begin{align*}
\frac{3 s}{1 + 3 s} \geq \frac{12}{5} s \geq 2 s \qquad \mbox{for all
  $s \in (0, 1/6)$. }
\end{align*}
Combining these bounds, we find that
\begin{align*}
  \P \Big[ \FDP_\numobs(t) \geq 2 s \quad \mbox{for all $t \in [0,
        \tau]$} \Big] & \geq \P \Big[\frac{L_\numobs\big ( \tau \big )
    }{L_\numobs\big ( \tau \big ) + n^{1 - \beta}} \geq \frac{3 s}{1 +
      3 s} \Big] \nonumber \\
& = \P \big[ L_\numobs(\tau) \geq 3 s \numobs^{1 - \beta} \big ].
\end{align*}
Consequently, the remainder of our proof is devoted to proving that
\begin{align}
\P \big[ L_\numobs(\tau) \geq 3 s \numobs^{1 - \beta} \big ] & \geq
1/2.
\end{align}
We split our analysis into two cases:

\paragraph{Case 1:}
First, suppose that $q_{n} \geq \frac{2 \log{4}}{3n^{1 - \beta}}$  In this case, we have
\begin{align}
  \label{eq:alpha-lower-bd}
\alpha \colon = \survival\big ( \tau \big ) \geq \frac{6 s}{n^{\beta}}
> \frac{16 \log{4}}{n}.
\end{align}
A simple calculation based on this inequality yields
\begin{align}
\alpha n - 3s\numobs^{1 - \beta} \geq \frac{\alpha n}{2}
\geq \sqrt{\big ( 4 \log 4 \big ) \alpha\big ( 1 - \alpha \big ) n} \defn a\sigma,
\end{align}
where $a = \sqrt{4 \log{4}}$ and $\sigma = \sqrt{\alpha\big ( 1 -
  \alpha \big ) n}$. Notice that $\sigma^{2} =
\Var\left[L_{\numobs}(\tau)\right]$.

We now apply the Bernstein inequality to $L_\numobs(\tau)$ to obtain
\begin{align*}
\P\big ( L_\numobs \leq 3s\numobs^{1 - \beta} \big ) & \leq \P\big (
L_\numobs \leq \pi_1 n - a\sigma \big ) \\
& \leq 2 \cdot \exp\big (
-\frac{a^{2}\sigma^{2}}{2\big[\sigma^{2} + a\sigma\big]} \big ) \\
& \leq 2 \cdot \exp\big ( -\frac{a^{2}}{2\big ( 1 + \frac{a}{\sigma}
  \big ) } \big ). \\ & \leq \exp\big ( -\frac{a^{2}}{4} \big ) =
\frac{1}{4},
\end{align*}
where we have used the fact that $a < \sigma$. We conclude that
\begin{align*}
\P\big ( L_\numobs \geq 3s\numobs^{1 - \beta} \big ) \geq \frac{1}{2},
\end{align*}
as desired.



\paragraph{Case 2:}  Otherwise, we may assume that
$q_{n} < \frac{2 \log{4}}{3n^{1 - \beta}}$. The definition of $\tau$ implies that
\begin{align*}
\alpha \geq \frac{24}{n} ~~ \text{and} ~~ 3s\numobs^{1 -
  \beta} \leq 8 \log{4}.
\end{align*}
It follows that $\E[L_\numobs(\tau)] \geq 24$. On the other hand,
given that $8 \log{4} < 12$, it suffices to prove that
\begin{align*}
P \big[ L_\numobs(\tau) \leq 12 \big] & \leq \frac{1}{2} .
\end{align*}
This is straightforward, however, since Bernstein's inequality gives
\begin{align*}
\P \big[ L_\numobs(\tau) \leq 12 \big] \; = \; \P \big[
  L_\numobs(\tau) \leq \frac{\E [L_n(\tau)]}{2} \big] & \leq \exp \big
( -\frac{24}{12} \big ) \\
& = e^{-2} \\
& < \frac{1}{2},
\end{align*}
which completes the proof.


\section{Proof of Lemmas~\ref{LemFNRFixed} and~\ref{lem:FNR-fixed-converse}}
\label{SecProofLemFNRFixed}

This appendix is devoted to the proofs of Lemmas~\ref{LemFNRFixed}
and~\ref{lem:FNR-fixed-converse}.  We combine the proofs, since these
two lemmas provide lower and upper bounds, respectively, on the
because they are matching lower and upper bounds, respectively, on the
FNR for a fixed threshold procedure, and their proofs involve
extremely similar calculations.

So as to simplify notation, we make use of the convenient shorthands
let $\tau = \taumin(q_\numobs)$, $\beta = \beta_\numobs$, and $\mu
= \mu_\numobs$ throughout the proof.  Recall that the FNP can be
written as the ratio $\FNP_\numobs(t) = \frac{F_\numobs(t)}{n^{1 -
    \beta}}$, where
\begin{align}
\label{eq:def-Fn}
F_\numobs(t) = \sum_{i \notin \nulls} \one\big ( X_{i} \leq t \big )
\sim \mathrm{Bin} \Big ( 1 - \survival\big ( t - \mu \big ) ,~n^{1 -
  \beta} \Big)
\end{align}
is a binomial random variable.  We split the remainder of the analysis
into two cases.

\paragraph{Case 1:}  First, suppose that $\tau \geq \mu$. In this case, we only seek to prove a lower bound. For this, observe that
$\survival\left(\tau - \mu\right) \leq \survival\left(0\right) = \frac{1}{2}$, so $1 - \survival\left(\tau - \mu\right) \geq \frac{1}{2}$. Thus,
\begin{align*}
\FDR_{n}\left(\tau\right) = \frac{\E\left[F_{\numobs}\right]}{n^{1 - \beta}} = 1 - \survival\left(\tau - \mu\right) \geq \frac{1}{2},
\end{align*}
as claimed.

\paragraph{Case 2:}  Otherwise, we may assume that $\mu > \tau$.
Recall the parameterization~\eqref{EqnMeanShift} of $\mu$ in terms of
$r$, the definition~\eqref{eq:rn-minus} of $\rmin$, and the
definition~\eqref{eq:dgamma-defn} of the $\dgammaexpr$ distance.  In
terms of these quantities, we have
\begin{align*}
\mu - \tau & = \big ( \gamma \log{n} \big )^{1/\gamma} \; \big \{
r^{1/\gamma} - \rmin(\kappan) \big \}^{1/\gamma} \\
& = \Big \{ \gamma \dgammaexpr \big(\rmin(\kappan), r \big) \;
\log{n} \Big \}^{1/\gamma} \\
& = \big[\gamma \dgammaexpr \Big(\beta + \kappan + \frac{\log
    \frac{1}{6\Zlower}}{\log n}, r \Big) \; \log{n}\big]^{1/\gamma},
\end{align*}
which shows how the quantity $\dgammaexpr$ determines the rate.  In
order to complete the proof, we need to show that the additional order
of $\frac{1}{\log n}$ term inside $\dgammaexpr$ can be removed.

More precisely, it suffices to establish the sandwich relation
\begin{align*}
\frac{\zeta^{2\beta^{\frac{1 - \gamma}{\gamma}}}}{\Zupper} \cdot
n^{-\dgamma{\beta + \kappan}{r}} \geq 1 - \survival\big ( \tau -
\mu \big ) \geq \frac{\zeta^{2\beta^{\frac{1 -
        \gamma}{\gamma}}}}{\Zlower} \cdot n^{-\dgamma{\beta +
    \kappan}{r}} ,
\end{align*}
where $\zeta = \max\big\lbrace 6\Zlower,~\frac{1}{6\Zlower}\big\rbrace$ as in~\eqref{EqZetaDef}. But now note that
\begin{align*}
\tau - \mu & = \big ( \gamma \log{n} \big ) ^{1/\gamma} \big[\big (
  \rmin \big ) ^{1/\gamma} - r\big] = -\big[\gamma
  \dgamma{\rmin}{r} \log{n}\big]^{1/\gamma},
\end{align*}
allowing us to deduce that
\begin{align*}
\frac{1}{\Zupper} \cdot n^{-\dgamma{\rmin}{r}} \geq 1 - \survival\big
( \tau - \mu \big ) \geq \frac{1}{\Zlower} \cdot
  n^{\dgamma{\rmin}{r}},
\end{align*}
so we need only show
\begin{align*}
  \big| \dgamma{\beta + \kappan}{r} - \dgamma{\rmin}{r}
  \big| \leq \frac{\beta^{\frac{1 - \gamma}{\gamma}} \log
    \zeta}{\log{n}} .
\end{align*}
To prove this, we let
\begin{align*}
\tilde{r}~\colon = \min\big ( \beta + \kappan,~\rmin \big ) 
\end{align*}
and note that by~\eqref{EqnAnnoying3}, we must have $\tilde{r} \in \left[\beta,~r\right]$. Under this definition, we consider the function $f(x) = \dgamma{\tilde{r} + x}{r}$. A
simple calculation shows that for $x \geq 0$, we have
\begin{align*}
f'(x) = \begin{cases} -\big ( \tilde{r} + x \big ) ^{\frac{1 -
      \gamma}{\gamma}}\dgamma{\tilde{r} + x}{r}^{\frac{\gamma -
      1}{\gamma}} &~ \text{if}~\tilde{r} + x \leq r, \\ \big (
  \tilde{r} + x \big ) ^{\frac{1 - \gamma}{\gamma}}\dgamma{\tilde{r} +
    x}{r}^{\frac{\gamma - 1}{\gamma}} &~\text{o.w.}
\end{cases}
\end{align*}
We observe that we only need to allow $0 \leq x \leq \max\big ( \beta
+ \kappan,~\rmin \big ) - \tilde{r} = \colon \tilde{R} -
\tilde{r}$, so in particular, we will always have $\tilde{r} + x \leq
\tilde{R} \leq 2$. This, together with the lower bound $\tilde{r} \geq
\beta$, yields
\begin{align*}
\sup_{0 \leq x \leq \tilde{R} - \tilde{r}}   \big|  f'(x) \big|  \leq 2\beta^{\frac{1 - \gamma}{\gamma}} .
\end{align*}
Applying this result, we find
\begin{align*}
  \big| \dgamma{\beta + \kappan}{r} - \dgamma{\rmin}{r}
  \big| & = \big| \dgamma{\tilde{R}}{r} - \dgamma{\tilde{r}}{r}
  \big| \\ & \leq 2\beta^{\frac{1 - \gamma}{\gamma}} \cdot \big (
  \tilde{R} - \tilde{r} \big ) \\ & = 2\beta^{\frac{1 -
      \gamma}{\gamma}} \cdot \frac{\log\zeta}{\log n} .
\end{align*}

If we now consider $q_{\numobs}' = cq_{\numobs}$, we can recover the more refined statements in Lemmas~\ref{LemFNRFixed} and~\ref{lem:FNR-fixed-converse}, simply by noting that
the same reasoning as above shows
\begin{align*}
  \big| \dgamma{\beta + \kappan}{r} - \dgamma{\beta + \kappan'}{r}
  \big| & \leq 2\beta^{\frac{1 - \gamma}{\gamma}} \cdot \frac{\big|\log
    c\big|}{\log n} ,
\end{align*}
concluding the argument.



\bibliographystyle{plainnat} \bibliography{paper}

\end{document}